\documentclass[12pt,oneside,reqno]{amsart}
\usepackage{txfonts}
\usepackage{bbm}
\usepackage{cases}
\usepackage{amsmath}
\usepackage{graphicx}
\usepackage{mathrsfs}
\usepackage{stmaryrd}
\usepackage{amsfonts}
\usepackage{enumerate,amsmath,amssymb,amsthm}
\numberwithin{equation}{section}

\newcommand{\be}{\begin{eqnarray}}
\newcommand{\ee}{\end{eqnarray}}
\newcommand{\ce}{\begin{eqnarray*}}
\newcommand{\de}{\end{eqnarray*}}
\newtheorem{theorem}{Theorem}[section]
\newtheorem{lemma}[theorem]{Lemma}
\newtheorem{remark}[theorem]{Remark}
\newtheorem{definition}[theorem]{Definition}
\newtheorem{proposition}[theorem]{Proposition}
\newtheorem{Examples}[theorem]{Example}
\newtheorem{corollary}[theorem]{Corollary}

\def\Re{{\mathrm{Re}}}
\def\Im{{\mathrm{Im}}}
\def\var{{\mathrm{var}}}

\def\eps{\varepsilon}

\def\p{\partial}

\def\[{{\Big[}}
\def\]{{\Big]}}
\def\<{{\langle}}
\def\>{{\rangle}}
\def\({{\Big(}}
\def\){{\Big)}}

\def\bx{{\mathbf{x}}}

\def\dif{{\mathord{{\rm d}}}}

\def\no{\nonumber}
\def\={&\!\!=\!\!&}
\def\bt{\begin{theorem}}
\def\et{\end{theorem}}
\def\bl{\begin{lemma}}
\def\el{\end{lemma}}
\def\br{\begin{remark}}
\def\er{\end{remark}}
\def\bd{\begin{definition}}
\def\ed{\end{definition}}
\def\bp{\begin{proposition}}
\def\ep{\end{proposition}}
\def\bc{\begin{corollary}}
\def\ec{\end{corollary}}
\def\bx{\begin{Examples}}
\def\ex{\end{Examples}}

\def\cJ{{\mathcal J}}

\def\cL{{\mathcal L}}
\def\cM{{\mathcal M}}

\def\cS{{\mathcal S}}
\def\cT{{\mathcal T}}

\def\mE{{\mathbb E}}

\def\mH{{\mathbb H}}
\def\mI{{\mathbb I}}

\def\mL{{\mathbb L}}
\def\mM{{\mathbb M}}
\def\mN{{\mathbb N}}

\def\mP{{\mathbb P}}

\def\mR{{\mathbb R}}
\def\mS{{\mathbb S}}

\def\mU{{\mathbb U}}

\def\mW{{\mathbb W}}

\def\bL{{\mathbf L}}

\def\sB{{\mathscr B}}

\def\sD{{\mathscr D}}

\def\sF{{\mathscr F}}
\def\sG{{\mathscr G}}

\def\sL{{\mathscr L}}
\def\sM{{\mathscr M}}

\def\geq{\geqslant}
\def\leq{\leqslant}

\allowdisplaybreaks

\begin{document}

\title{Uniqueness of stable-like processes}

\date{}

\date{}
\author{{Zhen-Qing Chen} \quad
and \quad {Xicheng Zhang}}

\address{Zhen-Qing Chen:
Department of Mathematics, University of Washington, Seattle, WA 98195, USA\\
Email: zqchen@uw.edu
 }
\address{Xicheng Zhang:
School of Mathematics and Statistics, Wuhan University,
Wuhan, Hubei 430072, P.R.China\\
Email: XichengZhang@gmail.com
 }

\begin{abstract}
In this work we consider the following $\alpha$-stable-like operator (a class of pseudo-differential operator)
$$
\sL f(x):=\int_{\mR^d}[f(x+\sigma_x y)-f(x)-1_{\alpha\in[1,2)}1_{|y|\leq 1}\sigma_x y\cdot\nabla f(x)]\nu_x(\dif y),
$$
where the L\'evy measure $\nu_x(\dif y)$ is comparable with a non-degenerate $\alpha$-stable-type L\'evy measure (possibly singular),
and $\sigma_x$ is a bounded and nondegenerate matrix-valued function. Under H\"older assumption on $x\mapsto\nu_x(\dif y)$ 
and uniformly continuity assumption on $x\mapsto\sigma_x$, we show the well-posedness of martingale problem associated with the operator $\sL$. Moreover, we also obtain the 
existence-uniqueness of strong solutions for the associated SDE when $\sigma$ belongs to the first order Sobolev space $\mW^{1,p}(\mR^d)$ provided
$p>d(1+\alpha\vee 1)$ and $\nu_x=\nu$ is  a non-degenerate $\alpha$-stable-type L\'evy measure.

\end{abstract}

\thanks{The research of ZC is partially supported
by NSF grant DMS-1206276. The research of XZ is partially supported by NNSFC grant of China (Nos. 11271294, 11325105).}

\maketitle

\section{Introduction}

Let $L_t$ be a $d$-dimensional L\'evy process. Consider the following stochastic differential equation (abbreviated as SDE) in $\mR^d$:
\begin{align}
\dif X_t=\sigma(X_{t-})\dif L_t,\ \ X_0=x,\label{SDE}
\end{align} 
where $\sigma:\mR^d\to\mR^d\otimes\mR^d$ is a measurable function. It is well-known that when $\sigma$ is Lipschitz continuous, 
there exists a unique strong solution $X_t(x)$ to SDE (\ref{SDE}) with infinitesimal generator
$$
\sL f(x):=\int_{\mR^d}[f(x+\sigma(x)y)-f(x)-1_{|y|\leq 1}\sigma(x)y\cdot\nabla f(x)]\nu(\dif y),
$$
where $\nu$ is the L\'evy measure of $L_t$.

\

Beyond the Lipschitz continuity assumption on $\sigma$, in the theory of SDEs, there are two different notions associated to the existence-uniqueness: 
weak existence-uniqueness (or uniqueness in law of weak solutions) and strong existence-uniqueness (or pathwise uniqueness of weak solutions).
Usually, strong existence-uniqueness requires stronger regularity conditions on $\sigma$ than weak existence-uniqueness.
When $L_t$ is a Brownian motion and $\sigma$ is uniformly non-degenerate and bounded {\it continuous}, in \cite{St-Va}
Stroock and Varadahan introduced the notion of martingale solutions, and studied the well-posedness of SDE (\ref{SDE}) in the weak sense.
In \cite{Kr}, Krylov established the existence of weak solutions to SDEs (\ref{SDE}) when $\sigma$ is uniformly nondegenerate and bounded {\it measurable}.
Recently, strong uniqueness for SDE (\ref{SDE}) driven by Brownian motion 
was proven in \cite{Zh0} when $\sigma$ is uniformly nondegenerate and belongs to the first order Sobolev space $\mW^{1,p}_{loc}$ provided $p>d$.

\

Nowadays, there has been a relatively complete theory for SDEs driven by Brownian motion. However, in the case of discontinuous L\'evy processes, 
there does not exist a satisfactory theory since the L\'evy measure of $L_t$ possess diversity and the associated Kolmogorov equation is {\it nonlocal}.
Such a feature brings us many difficulties. In particular, the study of the associated nonlocal integro-partial differential equation 
becomes more complicated, and strongly depends on the shape of the L\'evy measure. When $L_t$ is a cylindrical $\alpha$-stable process
and $\sigma$ is bounded continuous and nondegenerate, Bass and Chen \cite{Ba-Ch} proved the existence and uniqueness of weak solutions.
Notice that in this case, the L\'evy measure is singular with respect to the Lebesgue measure, and the argument is based on some singular-integral estimates of
pseudo-differential operators with singular state-dependent symbols.
When the L\'evy measure $\nu$ is absolutely continuous, or more generally,
$$
\sL' f(x):=\int_{\mR^d}[f(x+y)-f(x)-1_{|y|\leq 1}y\cdot\nabla f(x)]\frac{\kappa(x,y)}{|y|^{d+\alpha}}\dif y,
$$
where $\alpha\in(0,2)$, there are a lot of works devoting to the well-posedness of the martingale problem 
associated to $\sL'$ perturbed by some lower order term under different assumptions 
(see \cite{Ko, Mi-Pr1, Mi-Pr2, Ab-Ka} and references therein). To the best of our knowledge, the weakest assumptions on $\kappa(x,y)$ are given in \cite{Mi-Pr2}, 
i.e., for some $\kappa_0,\kappa_1>0$ and $\gamma\in(0,1)$,
$$
0<\kappa_0\leq\kappa(x,y)\leq\kappa_1,\ \ |\kappa(x,y)-\kappa(x',y)|\leq C|x-x'|^\gamma.
$$

The purpose of this work is to study the strong and weak uniqueness of SDEs associated to the following more general L\'evy operator:
\begin{align*}
\sL'' f(x)&:=\int_{\mR^d}[f(x+\sigma_xy)-f(x)-1_{|y|\leq 1}\sigma_xy\cdot\nabla f(x)]\nu_x(\dif y)\\
&+\int_{\mR^d}[f(x+\bar\sigma_xy)-f(x)-1_{|y|\leq 1}\bar\sigma_xy\cdot\nabla f(x)]\bar\nu_x(\dif y),
\end{align*}
where $\nu_x$ is comparable with a nondegenerate $\alpha$-stable L\'evy measure, 
and $\bar\nu_x$ is bounded by some $\beta$-stable L\'evy measure with $0<\beta<\alpha$, $\sigma_x$ is bounded continuous and non-degenerate, $\bar\sigma_x$
is bounded measurable. It should be observed that if $\nu_x(\dif y)$ and $\bar\nu_x(\dif y)$ are absolutely continuous with respect to the Lebesgue measure, 
then by the change of variables, the operator $\sL''$ can be written as $\sL'$ perturbed by some lower order term. 
Here we allow $\nu_x(\dif y)$ to be singular so that it can cover SDE (\ref{SDE}). Since the symbol of $\sL''$ could be very singular along the axis, 
we can not use the theory of the classical pseudo-differential operator to study the associated parabolic equation. 
We shall use the $L^p$-maximal regularity of nonlocal operator established recently in \cite{Zh1} to study the
solvability of the nonlocal parabolic equation associated with $\sL''$. After this, we shall study the existence and uniqueness of martingale solutions associated with $\sL''$
by establishing a Krylov's type estimate.

\

Compared with the existing literatures, the novelty of this work lies in the following two points:
\begin{enumerate}[(i)]
\item We do not assume that $\nu_x$ and $\bar\nu_x$ are absolutely continuous so that it can be used to solve the following SDE:
\begin{align}\label{SDE0}
\dif X_t=\sigma(X_t)\dif L_t+\bar\sigma(X_t)\dif \bar L_t,
\end{align}
where $L_t$ and $\bar L_t$ are two independent L\'evy processes, the L\'evy measure of $L_t$ is comparable with a nondegenerate $\alpha$-stable L\'evy measure,
and the L\'evy measure of $\bar L_t$ is just bounded by a $\beta$-stable L\'evy measure with $\beta<\alpha$. In particular, the cylindrical L\'evy processes are allowed.

\vspace{2mm}

\item We do not make any {\it H\"older} assumptions on $\sigma(x)$ and $\bar\sigma(x)$ so that our existence and uniqueness can cover most of the well-known results
such as the ones studied in \cite{Ba-Ch} and \cite{Mi-Pr2}.
\end{enumerate}

This paper is organised as follows: In Section 2, we give some preliminaries, and particularly, establish some auxiliary estimates. 
We believe that part of them has some independent interest (for example, Theorem \ref{Th3} below). 
In Section 3, we study a quite general nonlocal parabolic equation with space-time dependent coefficients,
and establish the $L^p$-maximal solvability by using Levi's freezing coefficient argument.
In Section 4, basing on the main result in Section 3, we prove a Krylov's type estimate for the martingale problem associated with the nonlocal operator 
and then obtain the well-posedness of the martingale problem.
Finally, in Section 5, we also prove a pathwise uniqueness result when $\sigma$ is nondegenerate and belongs to $\mW^{1,p}_{loc}(\mR^d)$ 
with $p>d(1+\alpha\wedge 1)$.

\

Convention: The letter $C$ with or without subscripts will denote an unimportant constant, whose value may change in different places. Moreover, $A\preceq B$
means that $A\leq CB$ for some constant $C>0$, and $A\asymp B$ means that  $C^{-1}B\leq A\leq CB$ for some $C>1$.

\section{Preliminaries}

In this section, we introduce necessary spaces and lemmas for later use.
\subsection{Sobolev spaces and embeddings}
For $\alpha\geq0$ and $p>1$, let $H^{\alpha,p}:=(I-\Delta)^{-\frac{\alpha}{2}}(L^p(\mR^d))$
be the usual Bessel potential space with the norm
$$
\|f\|_{\alpha,p}:=\|(I-\Delta)^{\frac{\alpha}{2}}f\|_p\asymp\|f\|_p+\|\Delta^{\frac{\alpha}{2}}f\|_p,
$$
where $\|\cdot\|_p$ denotes the usual $L^p$-norm and $\Delta^{\frac{\alpha}{2}}=-(-\Delta)^{\frac{\alpha}{2}}$ is the  fractional Laplacian. 
For $m\in\mN$, an equivalent norm in  $H^{m,p}$ is given by
$$
\|f\|_{m,p}:=\sum_{k=0}^{m}\|\nabla^k f\|_p,
$$
where $\nabla^k$ denotes the $k$-order gradient.
Notice that the following interpolation inequality holds (cf. \cite{Be-Lo}): for any $\beta\in(0,\alpha), p>1$ and $f\in H^{\alpha,p}$,
\begin{align}
\|\Delta^{\frac{\beta}{2}}f\|_p\leq C_{d,p,\alpha,\beta}\|f\|^{1-\frac{\beta}{\alpha}}_p
\|\Delta^{\frac{\alpha}{2}}f\|^{\frac{\beta}{\alpha}}_p\leq\eps\|\Delta^{\frac{\alpha}{2}}f\|_p+C_\eps\|f\|_p, \ \  \eps>0,\label{Inter}
\end{align}
and by the boundedness of Riesz's transformation in $L^p$-spaces (cf. \cite{Ste}),
\begin{align}
\|\nabla f\|_p\asymp\|\Delta^{\frac{1}{2}}f\|_p,\ \ p>1.\label{Riesz}
\end{align}

The following lemma can be found in \cite{Zh0} and \cite{Ste}.
\bl\label{Le21}
\begin{enumerate}[(i)]
\item There exists a constant $C=C(d)>0$ such that for any $f\in C^1(\mR^d)$ and all $x,y\in\mR^d$,
$$
|f(x)-f(y)|\leq C|x-y|(\cM|\nabla f|(x)+\cM|\nabla f|(y)),
$$
where $\cM|\nabla f|(x):=\sup_{r>0}\frac{1}{|D_r|}\int_{D_r}|\nabla f(x+z)|\dif z$, and $D_r:=\{x: |x|\leq r\}$.
\item For any $p>1$, the maximal operator $\cM$ is bounded from $L^p$ to $L^p$.
\end{enumerate}
\el
The following two embedding results are more or less well known. For the reader's convenience, we provide their proofs here.
\bl\label{Le27}
For any $m\in\mN$, $\beta\in(0,1)$ and $p>m/\beta$, there is a constant $C=C(m,p,\beta)>0$ such that for any $f:\mR^m\to L^p$ and
all $(y_0,\delta)\in\mR^m\times\mR_+$,
$$
\Big\|\sup_{|y-y_0|\leq\delta}|f(y,\cdot)-f(y_0,\cdot)|\Big\|_p\leq C\delta^{\beta}\sup_{y,y'\in D_\delta(y_0)}\frac{\|f(y,\cdot)-f(y',\cdot)\|_p}{|y-y'|^\beta},
$$
where $D_\delta(y_0):=\{y\in\mR^m: |y-y_0|\leq \delta\}$.
\el
\begin{proof}
By considering $y\mapsto f(y,\cdot)-f(y_0,\cdot)\in L^p$, we may assume $y_0=0$ and $f(y_0,\cdot)=0$.
Let $\gamma\in(\frac{m}{p},\beta)$.
By Garsia-Rademich-Rumsey's inequality (see \cite[Theorem 2.1.3]{St-Va} or \cite[Lemma 23.2]{Hu}), 
there is a constant $C=C(m,p,\gamma)>0$ such that for all $x\in\mR^d$,
$$
\sup_{|y|\leq\delta}|f(y,x)|^p\leq C\delta^{\gamma p-m}\int_{D_\delta(0)}\!\int_{D_\delta(0)} \frac{|f(y,x)-f(y',x)|^p}{|y-y'|^{m+\gamma p}}\dif y\dif y'.
$$
Integrating both sides with respect to $x$, we get
\begin{align*}
\Big\|\sup_{|y|\leq\delta}|f(y,\cdot)\Big\|^p_p&\preceq\delta^{\gamma p-m}\int_{D_\delta(0)}\!\int_{D_\delta(0)} \frac{\|f(y,\cdot)-f(y',\cdot)\|^p_p}{|y-y'|^{m+\gamma p}}\dif y\dif y'\\
&\preceq K^p\delta^{\gamma p-m}\int_{D_\delta(0)}\!\int_{D_\delta(0)} |y-y'|^{-m+(\beta-\gamma) p}\dif y\dif y'
\preceq K^p\delta^{\beta p},
\end{align*}
where $K:=\sup_{y,y'\in D_\delta(0)}\frac{\|f(y,\cdot)-f(y',\cdot)\|_p}{|y-y'|^\beta}$. The proof is complete.
\end{proof}
Below, for $0<S<T$, we shall write 
$$
\mL^p(S,T):=L^p([S,T]\times\mR^d)=L^p([S,T]; L^p(\mR^d))
$$
and
$$
\mH^{\alpha,p}(S,T):=L^p([S,T]; H^{\alpha,p}).
$$
\bl\label{Le23}
For any $p>1$ and $\beta\in(0,\alpha(1-\frac{1}{p}))$, there exists a constant $C=C(d,p,\alpha,\beta)>0$ 
such that for all $t_0<{t_1}$,
$$
\|u({t_1})-u(t_0)\|_{ \beta,p}\leq C({t_1}-t_0)^{1-\frac{\beta}{\alpha}-\frac{1}{p}}\Big(\|\p_t u\|_{\mL^p(t_0,{t_1})}+\|u\|_{\mH^{\alpha,p}(t_0,{t_1})}\Big),
$$
provided that the right hand side is finite.
\el
\begin{proof}
Since $(1-\frac{\beta}{\alpha})p>1$, one can choose 
$$
\gamma\in(0, 1-\tfrac{\beta}{\alpha}),\ \ \delta\in(1-\tfrac{(\alpha(1-\gamma)-\beta)p}{\beta},1)
$$ 
such that
\begin{align}\label{ES4}
\gamma p>1,\ \ (\alpha-\beta)(p+1)>\alpha+\alpha\gamma p-\delta\beta.
\end{align}
By Garsia-Rademich-Rumsey's inequality again, there exits a constant $C=C(\gamma,p)>0$ such that for all $t_0<{t_1}$,
\begin{align}\label{ES2}
\|u({t_1})-u(t_0)\|^p_{\beta,p}\leq C({t_1}-t_0)^{\gamma p-1}\int^{t_1}_{t_0}\!\!\int^t_{t_0}\frac{\|u(t)-u(s)\|^p_{ \beta,p}}{(t-s)^{1+\gamma p}}\dif s\dif t.
\end{align}
By the interpolation inequality \eqref{Inter} and H\"older's inequality, we have
\begin{align*}
&\int^{t_1}_{t_0}\!\!\!\int^t_{t_0}\frac{\|u(t)-u(s)\|^p_{ \beta,p}}{(t-s)^{1+\gamma p}}\dif s\dif t\\
&\preceq\int^{t_1}_{t_0}\!\!\!\int^t_{t_0}\frac{\|u(t)-u(s)\|^{(\alpha-\beta)p/\alpha}_{p}\|u(t)-u(s)\|^{\beta p/\alpha}_{\alpha,p}}{(t-s)^{1+\gamma p}}\dif s\dif t\\
&\leq\left(\int^{t_1}_{t_0}\!\!\!\int^t_{t_0}\frac{\|u(t)-u(s)\|^{p}_{p}}{(t-s)^{((1+\gamma p)\alpha-\delta\beta)/(\alpha-\beta)}}\dif s\dif t\right)^{(\alpha-\beta)/\alpha}\\
&\quad\times\left(\int^{t_1}_{t_0}\!\!\!\int^t_{t_0}\frac{\|u(t)-u(s)\|^{p}_{\alpha,p}}{(t-s)^{\delta}}\dif s\dif t\right)^{\beta/\alpha}=:I_1\times I_2.
\end{align*}
To treat $I_1$, we need the following elementary estimate: for any $q>0$,
\begin{align}\label{ES1}
\int^{t_1}_{t_0}\!\!\!\int^t_{t_0}\left((t-s)^{q-2}\int^t_s f(r)\dif r\right)\dif s\dif t\leq \frac{({t_1}-{t_0})^q}{q|1-q|}\|f\|_{L^1({t_0},{t_1})}.
\end{align}
Indeed, let
$$
\sD:=\left\{t\in({t_0},{t_1}): \lim_{s\uparrow t}\frac{(t-s)^{q}}{t-s}\int^t_s f(r)\dif r=0\right\}.
$$
Since $f\in L^1({t_0},{t_1})$, by the Lebesgue differential theorem, $\sD$ has full measure. Thus, 
for each $t\in\sD$, by the integration by parts formula, we have
$$
\int^t_{t_0}\left((t-s)^{q-2}\int^t_s f(r)\dif r\right)\dif s=\frac{1}{|1-q|}\int^t_{t_0}(t-s)^{q-1}f(s)\dif s.
$$
Hence,
$$
\int^{t_1}_{t_0}\!\!\!\int^t_{t_0}\left((t-s)^{q-2}\int^t_s f(r)\dif r\right)\dif s\dif t=\frac{1}{|1-q|}\int^{t_1}_{t_0}\!\!\!\int^t_{t_0}(t-s)^{q-1}f(s)\dif s\dif t,
$$
which in turn implies \eqref{ES1} by Fubini's theorem. 

Now, noticing that
$$
\|u(t)-u(s)\|^{p}_{p}\leq (t-s)^{p-1}\int^t_s\|\p_r u(r)\|^{p}_{p}\dif r,
$$
by \eqref{ES1} and \eqref{ES4}, we have
$$
I_1\preceq ({t_1}-{t_0})^{(1-\frac{\beta}{\alpha}-\gamma) p-\frac{(1-\delta)\beta}{\alpha}}\|\p_t u\|_{\mL^p({t_0},{t_1})}^{(\alpha-\beta)p/\alpha}.
$$
For $I_2$, we have
\begin{align*}
I_2&\preceq\left(\int^{t_1}_{t_0}\|u(t)\|^p_{\alpha,p}\!\!\int^t_{t_0}\frac{\dif s}{(t-s)^{\delta}}\dif t
+\int^{t_1}_{t_0}\|u(s)\|^p_{\alpha,p}\!\!\int^{t_1}_s\frac{\dif t}{(t-s)^{\delta}}\dif s\right)^{\beta/\alpha}\\
&\preceq({t_1}-{t_0})^{\frac{(1-\delta)\beta}{\alpha}}\left(\int^{t_1}_{t_0}\|u(t)\|^p_{\alpha,p}\dif t\right)^{\beta/\alpha}
=({t_1}-{t_0})^{\frac{(1-\delta)\beta}{\alpha}}\|u\|_{ \mH^{\alpha,p}({t_0},{t_1})}^{\beta p/\alpha}.
\end{align*}
Hence,
$$
\int^{t_1}_{t_0}\!\!\!\int^t_{t_0}\frac{\|u(t)-u(s)\|^p_{ \beta,p}}{(t-s)^{1+\gamma p}}\dif s\dif t\preceq ({t_1}-{t_0})^{(1-\frac{\beta}{\alpha}-\gamma)p}
\|\p_t u\|_{\mL^p({t_0},{t_1})}^{(\alpha-\beta)p/\alpha}\|u\|_{\mH^{\alpha,p}({t_0},{t_1})}^{\beta p/\alpha},
$$
which together with \eqref{ES2} gives the desired estimate.
\end{proof}

For $\alpha\in(0,2)$ and $y\in\mR^d$,  we write
$$
y^{(\alpha)}:=y1_{|y|\leq 1}1_{\alpha=1}+y1_{\alpha\in(1,2)},
$$
and for a function $f:\mR^d\to\mR$,
\begin{align}
\cJ^{(\alpha)}_f(x,y):=f(x+y)-f(x)-y^{(\alpha)}\cdot\nabla f(x).\label{JJ}
\end{align}
The following lemma is taken from \cite[Lemma 5]{Mi-Pr1}.
\bl\label{Le22}
For $\alpha\in(0,2)$ and $p>\frac{d}{\alpha}\vee 1$, there is a constant $C=C(p,d,\alpha)>0$ such that for all $f\in H^{\alpha,p}$,
\begin{align}
\left\|\sup_{y\not=0}\frac{|\cJ^{(\alpha)}_f(\cdot,y)|}{|y|^\alpha}\right\|_p\leq C\|\Delta^{\frac{\alpha}{2}}f\|_p.\label{Es34}
\end{align}
\el
The following lemma is direct by Sobolev's embedding theorem.
\bl
For $\alpha\in(0,2)$, $\beta\in(\alpha,2)$ and $p>\frac{d}{\beta-\alpha}\vee 1$, there is a constant $C=C(p,d,\alpha,\beta)>0$ such that for all $f\in H^{\beta,p}$,
\begin{align}
\sup_x\sup_{y\not=0}\frac{|\cJ^{(\alpha)}_f(x,y)|}{|y|^\alpha}\leq C\|f\|_{\beta,p}.\label{Es304}
\end{align}
\el
\subsection{$L^p$-estimate of L\'evy operators}
Let $\bL$ be the set of all L\'evy measures $\nu$ on $\mR^d$, that is,
$$
\nu(\{0\})=0,\ \ \int_{\mR^d}1\wedge|x|^2\nu(\dif x)<+\infty,
$$
which is endowed with the weak convergence topology.
For $\alpha\in(0,2)$, let $\bL^{(\alpha)}\subset\bL$ be the set of all $\alpha$-stable measure $\nu^{(\alpha)}$ with the form
\begin{align}
\nu^{(\alpha)}(\Gamma):=\int^\infty_0\left(\int_{\mS^{d-1}}\frac{1_\Gamma(r\theta)\Sigma(\dif\theta)}{r^{1+\alpha}}\right)\dif r,\ \ \Gamma\in\sB(\mR^d),\label{Eq4}
\end{align}
where $\Sigma$ is a finite measure over the sphere $\mS^{d-1}$ (called spherical measure of $\nu^{(\alpha)}$), and we also require that
$$
1_{\alpha=1}\int_{\mS^{d-1}}\theta\Sigma(\dif\theta)=0.
$$
Let $\mM^d$ be the space of all real invertible $d\times d$-matrix. The identity matrix is denoted by $\mI$, and the transpose of a matrix $\sigma$ is denoted by $\sigma^*$.
Let $\cS(\mR^d)$  be the Schwartz rapidly decreasing function space. 

Given $\nu\in\bL$, $\sigma\in\mM^d$ and $\alpha\in(0,2)$, we consider the following L\'evy operator:
$$
\cL^\nu_\sigma f(x):=\int_{\mR^d}\cJ^{(\alpha)}_f(x,\sigma y)\nu(\dif y),\ \ f\in\cS(\mR^d),
$$
where $\cJ^{(\alpha)}_f(x,\sigma y)$ is defined by (\ref{JJ}). 
Clearly,
\begin{align}\label{For}
\begin{split}
&\cL^\nu_\sigma (fg)(x)-f(x)\cL^\nu_\sigma g(x)-g(x)\cL^\nu_\sigma f(x)\\
&\quad=\int_{\mR^d}(f(x+\sigma y)-f(x))(g(x+\sigma y)-g(x))\nu(\dif y).
\end{split}
\end{align}
Let $\psi^\nu_\sigma$ be the symbol of operator $\cL^\nu_\sigma$, i.e.,
$$
\widehat{\cL^\nu_\sigma f}(\xi)=\psi^\nu_\sigma(\xi) \hat f(\xi),
$$
where $\hat{f}$ denotes Fourier's transformation of $f$. It is easy to see that
\begin{align}
\psi^\nu_\sigma(\xi):=\int_{\mR^d}(1+\mathrm{i}\xi\cdot \sigma y^{(\alpha)}-\mathrm{e}^{\mathrm{i}\xi\cdot \sigma y})\nu(\dif y).\label{psi}
\end{align}
In particular, if $\nu(\dif y)=|y|^{-d-\alpha}\dif y$, then
\begin{align}
\psi^\nu_\mI(\xi)=c_{d,\alpha}|\xi|^\alpha\ \mbox{ and }\ \cL^\nu_\mI f(x)=c_{d,\alpha}\Delta^{\frac{\alpha}{2}}f(x).\label{Es7}
\end{align}

We introduce the following notions.
\bd
\begin{enumerate}[(i)]
\item For $\nu^{(\alpha)}\in \bL^{(\alpha)}$, it is called non-degenerate if
$$
\int_{\mS^{d-1}}|\theta_0\cdot\theta|^\alpha\Sigma(\dif\theta)\not=0,\ \ \forall\theta_0\in\mS^{d-1},
$$
where $\Sigma$ is the spherical measure of $\nu^{(\alpha)}$. The set of all non-degenerate $\alpha$-stable measures is denoted by $\bL^{(\alpha)}_{non}$.
\item For $\nu_1,\nu_2\in\bL$, we say that $\nu_1$ is less than $\nu_2$ if
$$
\nu_1(\Gamma)\leq \nu_2(\Gamma),\ \ \Gamma\in\sB(\mR^d),
$$
and we simply write $\nu_1\leq \nu_2$ in this case.
\end{enumerate}
\ed
The following lemma gives a characterization of non-degenerate L\'evy measures.
\bl\label{Le11}
Let $\nu\in\bL$, $\nu^{(\alpha)}\in\bL^{(\alpha)}_{non}$ and $\sigma\in\mM^d$. If $\nu\geq \nu^{(\alpha)}$, then
\begin{align}
\mathrm{Re}(\psi^\nu_\sigma(\xi))\geq c_\alpha\left(\inf_{\theta\in\mS^{d-1}}|\sigma^*\theta|^\alpha\inf_{\theta_0\in\mS^{d-1}}
\int_{\mS^{d-1}}|\theta_0\cdot \theta|^\alpha\Sigma(\dif\theta)\right)|\xi|^\alpha,\label{EP8}
\end{align}
where $c_\alpha$ only depends on $\alpha$, and $\Sigma$ is the spherical measure of $\nu^{(\alpha)}$.
\el
\begin{proof}
By (\ref{psi}) and the change of variables, we have
\begin{align*}
\mathrm{Re}(\psi^\nu_\sigma(\xi))&=\int_{\mR^d}(1-\cos(\xi\cdot \sigma y))\nu(\dif y)\geq\int_{\mR^d}(1-\cos(\xi\cdot \sigma y))\nu^{(\alpha)}(\dif y)\\
&=\int^\infty_0\left(\int_{\mS^{d-1}}\frac{(1-\cos(r\xi\cdot \sigma \theta))\Sigma(\dif\theta)}{r^{1+\alpha}}\right)\dif r\\
&=\left(\int^\infty_0\frac{1-\cos r}{r^{1+\alpha}}\dif r\right)\left(\int_{\mS^{d-1}}|\xi\cdot \sigma \theta|^\alpha\Sigma(\dif\theta)\right)\\
&\geq\left(\int^\infty_0\frac{1-\cos r}{r^{1+\alpha}}\dif r\right)\left(\inf_{\theta_0\in\mS^{d-1}}\int_{\mS^{d-1}}|\theta_0\cdot \theta|^\alpha\Sigma(\dif\theta)\right)|\sigma^*\xi|^\alpha,
\end{align*}
which then gives (\ref{EP8}).
\end{proof}
Next we show the continuous dependence of the symbol $\psi^\nu_\sigma$ with respect to $\nu$ and $\sigma$. 
We need the following elementary estimate.
\bl\label{Le5}
Let $a,b\in\mR$. We have
\begin{enumerate}[(i)]
\item If $\alpha\in(0,1)$, then
$$
\int^\infty_0\Big(|\cos (ar)-\cos(br)|+|\sin(ar)-\sin(br)|\Big)\frac{\dif r}{r^{1+\alpha}}\leq c_\alpha |a-b|^\alpha.
$$
\item If $\alpha=1$, then for any $\beta\in(0,1)$,
\begin{align*}
\int^\infty_0\Big(|\cos (ar)-\cos(br)|&+|(a-b)r1_{r\leq |a-b|^{-1}}-(\sin(ar)-\sin(br))|\Big)\frac{\dif r}{r^{2}}\\
&\leq c_\beta (|a|+|b|)^{1-\beta}~|a-b|^{\beta},
\end{align*}
\item If $\alpha\in(1,2)$, then
\begin{align*}
\int^\infty_0\Big(|\cos (ar)-\cos(br)|&+|(a-b) r-(\sin(ar)-\sin(br))|\Big)\frac{\dif r}{r^{1+\alpha}}\\
&\leq c_\alpha (|a|+|b|)^{\alpha-1}~|a-b|,
\end{align*}
\end{enumerate}
Here $c_\alpha$ and $c_\beta$ only depends on $\alpha$ and $\beta$.
\el
\begin{proof}
Below, we assume $a\not=b$.
\\
(i) By $|\cos x-\cos y|\leq|x-y|$ and $|\sin x-\sin y|\leq|x-y|$, we have
\begin{align*}
&\int^\infty_0(|\cos (ar)-\cos(br)|+|\sin(ar)-\sin(br)|)\frac{\dif r}{r^{1+\alpha}}\\
&\quad\leq 2\int^{|a-b|^{-1}}_0|a-b|\frac{\dif r}{r^{\alpha}}+4\int^\infty_{|a-b|^{-1}}\frac{\dif r}{r^{1+\alpha}}\leq c_\alpha|a-b|^\alpha.
\end{align*}
(ii) By $|\sin x|\leq |x|^\beta$, we have
\begin{align*}
&\int^\infty_0|\cos (ar)-\cos(br)|\frac{\dif r}{r^{2}}\leq 2\int^\infty_{|a-b|^{-1}}\frac{\dif r}{r^2}\\
&\quad+|a-b|\int^{|a-b|^{-1}}_0\!\!\!\!\!\int^1_0|\sin(r(a(1-s)+bs))|\dif s\frac{\dif r}{r}\\
&\quad\leq2|a-b|+|a-b|(|a|+|b|)^{1-\beta}\int^{|a-b|^{-1}}_0 r^{-\beta}\dif r\\
&\quad\leq c_\beta|a-b|^{\beta}(|a|+|b|)^{1-\beta},
\end{align*}
and by $|1-\cos x|\leq |x|^\beta$,
\begin{align*}
&\int^\infty_0|(a-b)r1_{r\leq |a-b|^{-1}}-(\sin(ar)-\sin(br))|\frac{\dif r}{r^{2}}\\
&\quad\leq |a-b|\int^{|a-b|^{-1}}_0\!\!\!\int^1_01-\cos(r(a(1-s)+bs))\dif s\frac{\dif r}{r}+2\int^\infty_{|a-b|^{-1}}\frac{\dif r}{r^2}\\
&\quad\leq |a-b|(|a|+|b|)^{1-\beta}\int^{|a-b|^{-1}}_0r^{-\beta}\dif r+2|a-b|\leq c_\beta|a-b|^{\beta}(|a|+|b|)^{1-\beta}.
\end{align*}
(iii) By $|\sin x|\leq|x|$, we have
\begin{align*}
&\int^\infty_0\frac{|\cos (ar)-\cos(br)|}{r^{1+\alpha}}\dif r\\
&\leq|a-b|\int^\infty_0\!\!\!\int^1_0|\sin(r(a(1-s)+bs))|\dif s\frac{\dif r}{r^\alpha}\\
&\leq|a-b|\left((|a|+|b|)\int^{(|a|+|b|)^{-1}}_0\frac{r\dif r}{r^\alpha}+\int^\infty_{(|a|+|b|)^{-1}}\frac{\dif r}{r^\alpha}\right)\\
&\leq c_\alpha|a-b|(|a|+|b|)^{\alpha-1},
\end{align*}
and by $|1-\cos x|\leq|x|$,
\begin{align*}
&\int^\infty_0|(a-b) r-(\sin(ar)-\sin(br))|\frac{\dif r}{r^{1+\alpha}}\\
&\quad\leq|a-b|\int^\infty_0\!\!\!\int^1_0(1-\cos(r(a(1-s)+bs)))\frac{\dif r}{r^{\alpha}}\\
&\quad\leq|a-b|\left((|a|+|b|)\int^{(|a|+|b|)^{-1}}_0\frac{r\dif r}{r^\alpha}+2\int^\infty_{(|a|+|b|)^{-1}}\frac{\dif r}{r^\alpha}\right)\\
&\quad\leq c_\alpha|a-b|(|a|+|b|)^{\alpha-1}.
\end{align*}
The proof is complete.
\end{proof}

Now we can show the following continuous dependence of symbol $\psi^\nu_\sigma$ with respect to $\nu$ and $\sigma$.
\bl\label{Le12}
Let $\sigma_1,\sigma_2\in\mM^d$ and $\nu_1,\nu_2\in\bL$. Assume that for some $\nu^{(\alpha)}\in\bL^{(\alpha)}$ and $K>0$,
\begin{align}
\nu_1,\nu_2\leq\nu^{(\alpha)},\ \ |\nu_1-\nu_2|\leq K\nu^{(\alpha)},\label{NU01}
\end{align}
and for all $0<r<R<\infty$,
\begin{align}
1_{\alpha=1}\int_{r<|y|<R}y\nu_i(\dif y)=0,\ \ i=1,2.\label{NU1}
\end{align}
\begin{enumerate}[(i)]
\item
If $\alpha\in(0,1)$, then
$$
|\psi^{\nu_1}_{\sigma_1}(\xi)-\psi^{\nu_2}_{\sigma_2}(\xi)|\leq C(K+|\sigma_1-\sigma_2|^\alpha)|\xi|^\alpha.
$$
\item If $\alpha=1$, then for any $\beta\in(0,1)$,
$$
|\psi^{\nu_1}_{\sigma_1}(\xi)-\psi^{\nu_2}_{\sigma_2}(\xi)|\leq C(K+|\sigma_1-\sigma_2|^{\beta}(|\sigma_1|+|\sigma_2|)^{1-\beta})|\xi|^\alpha.
$$
\item If $\alpha\in(1,2)$, then
$$
|\psi^{\nu_1}_{\sigma_1}(\xi)-\psi^{\nu_2}_{\sigma_2}(\xi)|\leq C(K+|\sigma_1-\sigma_2|(|\sigma_1|+|\sigma_2|)^{\alpha-1})|\xi|^\alpha.
$$
\end{enumerate}
Here the constant $C$ only depends on $d,\alpha,\nu^{(\alpha)}$ and $\beta$.
\el
\begin{proof}
By (\ref{psi}), (\ref{NU01}), (\ref{Eq4}) and the change of variables, we have
\begin{align*}
|\Re(\psi^{\nu_1}_{\sigma_1}(\xi)-\psi^{\nu_2}_{\sigma_1}(\xi))|&\leq\int_{\mR^d}(1-\cos(\xi\cdot \sigma_1 y))\|\nu_1-\nu_2\|_\var(\dif y)\\
&\leq K\int_{\mR^d}(1-\cos(\xi\cdot \sigma_1 y))\nu^{(\alpha)}(\dif y)\\
&\leq K|\sigma_1|^\alpha|\xi|^\alpha\left(\int^\infty_0\frac{1-\cos r}{r^{1+\alpha}}\dif r\right)\Sigma(\mS^{d-1}),
\end{align*}
and
\begin{align*}
&|\Re(\psi^{\nu_2}_{\sigma_1}(\xi)-\psi^{\nu_2}_{\sigma_2}(\xi))|\\
&\leq\int_{\mR^d}|\cos(\xi\cdot \sigma_1 y)-\cos(\xi\cdot \sigma_2 y)|\nu^{(\alpha)}(\dif y)\\
&=\int_{\mS^{d-1}}\left(\int^\infty_0\frac{|\cos(r\xi\cdot \sigma_1\theta)-\cos(r\xi\cdot \sigma_2\theta)|\dif r}{r^{1+\alpha}}\right)\Sigma(\dif\theta),
\end{align*}
which implies the desired estimate for the real part of $\psi^{\nu_1}_{\sigma_1}(\xi)-\psi^{\nu_2}_{\sigma_2}(\xi)$ by Lemma \ref{Le5}.

On the other hand, if $\alpha\in(0,1)$, then
\begin{align*}
|\Im(\psi^{\nu_1}_{\sigma_1}(\xi)-\psi^{\nu_2}_{\sigma_1}(\xi))|&\leq K\int_{\mR^d}|\sin(\xi\cdot\sigma_1 y)|\nu^{(\alpha)}(\dif y)\\
&\leq K|\sigma_1|^\alpha|\xi|^\alpha\int_{\mS^{d-1}}\left(\int^\infty_0\frac{|\sin(r\widehat{\sigma^*_1\xi}\cdot  \theta)|\dif r}{r^{1+\alpha}}\right)\Sigma(\dif\theta)\\
&\leq K|\sigma_1|^\alpha|\xi|^\alpha\left(\int^\infty_0\frac{(r\wedge 1)\dif r}{r^{1+\alpha}}\right)\Sigma(\mS^{d-1}),
\end{align*}
where $\widehat{\sigma^*_1\xi}=\sigma^*_1\xi/|\sigma^*_1\xi|$, and
\begin{align*}
|\Im(\psi^{\nu_2}_{\sigma_1}(\xi)-\psi^{\nu_2}_{\sigma_2}(\xi))|&\leq\int_{\mR^d}|\sin(\xi\cdot\sigma_1 y)-\sin(\xi\cdot \sigma_2 y)|\nu^{(\alpha)}(\dif y)\\
&\leq\int_{\mS^{d-1}}\left(\int^\infty_0\frac{|r \xi\cdot\sigma_1 \theta-\sin(r\xi\cdot \sigma_2 \theta)|\dif r}{r^{1+\alpha}}\right)\Sigma(\dif\theta);
\end{align*}
if $\alpha=1$, then by (\ref{NU1}),
\begin{align*}
&|\Im(\psi^{\nu_1}_{\sigma_1}(\xi)-\psi^{\nu_2}_{\sigma_1}(\xi))|\\
&=\left|\int_{\mR^d}(\xi\cdot\sigma_1y 1_{|y|\leq|\sigma^*_1\xi|^{-1}}-\sin(\xi\cdot\sigma_1 y))(\nu_1-\nu_2)(\dif y)\right|\\
&\leq K\int_{\mR^d}|\xi\cdot\sigma_1y 1_{|y|\leq|\sigma^*_1\xi|^{-1}}-\sin(\xi\cdot\sigma_1 y))|\nu_1^{(\alpha)}(\dif y)\\
&\leq K|\sigma_1|^\alpha|\xi|^\alpha\int_{\mS^{d-1}}\left(\int^\infty_0\frac{|r\widehat{\sigma^*_1\xi}\cdot\theta 1_{r\leq 1}-\sin(r\widehat{\sigma^*_1\xi}\cdot \theta)|\dif r}{r^{1+\alpha}}\right)\Sigma(\dif\theta),
\end{align*}
and
\begin{align*}
&|\Im(\psi^{\nu_2}_{\sigma_1}(\xi)-\psi^{\nu_2}_{\sigma_2}(\xi))|\\
&\quad\leq\int_{\mR^d}|(\xi\cdot(\sigma_1-\sigma_2) y1_{|y|\leq |(\sigma_1-\sigma_2)^*\xi|^{-1}}\\
&\qquad\qquad-(\sin(\xi\cdot \sigma_1 y)-\sin(\xi\cdot \sigma_2 y))|\nu^{(\alpha)}(\dif y)\\
&\quad=\int_{\mS^{d-1}}\Bigg(\int^\infty_0|(r\xi\cdot(\sigma_1-\sigma_2) \theta1_{r\leq |(\sigma_1-\sigma_2)^*\xi|^{-1}}\\
&\qquad\qquad-(\sin(r\xi\cdot \sigma_1\theta)-\sin(r\xi\cdot \sigma_2\theta))|\frac{\dif r}{r^{1+\alpha}}\Bigg)\Sigma(\dif\theta);
\end{align*}
if $\alpha\in(1,2)$, then
\begin{align*}
&|\Im(\psi^{\nu_1}_{\sigma_1}(\xi)-\psi^{\nu_2}_{\sigma_1}(\xi))|\\
&\leq K\int_{\mR^d}|\xi\cdot\sigma_1y-\sin(\xi\cdot\sigma_1 y)|\nu^{(\alpha)}(\dif y)\\
&\leq K|\sigma_1|^\alpha|\xi|^\alpha\int_{\mS^{d-1}}\left(\int^\infty_0\frac{|r\widehat{\sigma^*_1\xi}\cdot\theta-\sin(r\widehat{\sigma^*_1\xi}\cdot \theta)|\dif r}{r^{1+\alpha}}\right)\Sigma(\dif\theta),
\end{align*}
and
\begin{align*}
&|\Im(\psi^{\nu_2}_{\sigma_1}(\xi)-\psi^{\nu_2}_{\sigma_2}(\xi))|\\
&\leq\int_{\mR^d}|(\xi\cdot\sigma_1 y-\sin(\xi\cdot \sigma_1 y))-(\xi\cdot\sigma_2 y-\sin(\xi\cdot \sigma_2 y))|\nu^{(\alpha)}(\dif y)\\
&=\int^\infty_0\left(\int_{\mS^{d-1}}\frac{|\cos(r\xi\cdot \sigma_1\theta)-\cos(r\xi\cdot \sigma_2\theta)|\Sigma(\dif\theta)}{r^{1+\alpha}}\right)\dif r.
\end{align*}
Combining the above calculations, and by Lemma \ref{Le5} again, we obtain the desired estimate for the image part of $\psi^{\nu_1}_{\sigma_1}(\xi)-\psi^{\nu_2}_{\sigma_2}(\xi)$.
\end{proof}

Using Lemmas \ref{Le11} and \ref{Le12}, the following results are proven in \cite[Theorem 4.3 and Corollary 4.4]{Zh1}.
\bt\label{Th3}
\begin{enumerate}[(i)]
\item Let $\nu\in\bL$, $\nu^{(\alpha)}_1,\nu^{(\alpha)}_2\in \bL^{(\alpha)}_{non}$ and $\sigma\in\mM^d$, $\kappa_1,\kappa_2>0$.  Assume that
$$
\nu^{(\alpha)}_1\leq\nu\leq\nu^{(\alpha)}_2,\ \ 1_{\alpha=1}\int_{r<|y|<R}y\nu(\dif y)=0,\ \ 0<r<R<\infty,
$$
and
$$
\kappa_1|\xi|\leq |\sigma^*\xi|\leq \kappa_2|\xi|,\ \ \xi\in\mR^d.
$$
For any $p>1$, there exists a constant $C_0>0$ only depending on $d,p,\alpha,\nu^{(\alpha)}_1,\nu^{(\alpha)}_2,\kappa_1,\kappa_2$ such that for all $f\in H^{\alpha,p}$,
\begin{align}
C_0\|\Delta^{\frac{\alpha}{2}}f\|_p\leq\|\cL^\nu_\sigma f\|_p\leq C^{-1}_0\|\Delta^{\frac{\alpha}{2}}f\|_p.\label{EW888}
\end{align}
\item Let $\sigma_1,\sigma_2\in\mM^d$ and $\nu_1,\nu_2\in\bL$. Assume that for some $\nu^{(\alpha)}\in\bL^{(\alpha)}$ and $K, \kappa>0$,
\begin{align}
\nu_1,\nu_2\leq\nu^{(\alpha)},\ \ |\nu_1-\nu_2|\leq K\nu^{(\alpha)},\ \ |\sigma_1|,|\sigma_2|\leq\kappa,\label{NU}
\end{align}
and for all $0<r<R<\infty$,
\begin{align}
1_{\alpha=1}\int_{r<|y|<R}y\nu_i(\dif y)=0,\ \ i=1,2.\label{NU2}
\end{align}
For any $p>1$ and $\beta\in(0,1)$, there is a constant $C_1>0$ only depending on $d,p,\alpha,\nu^{(\alpha)},\kappa,\beta$ such that for any $f\in H^{\alpha,p}$,
\begin{align}\label{EW88}
\|\cL^{\nu_1}_{\sigma_1} f-\cL^{\nu_2}_{\sigma_2} f\|_p\leq C_1(K+|\sigma_1-\sigma_2|^{\beta_\alpha})\|\Delta^{\frac{\alpha}{2}}f\|_p,
\end{align}
where 
\begin{align}
\beta_\alpha:=\alpha 1_{\alpha\in(0,1)}+\beta 1_{\alpha=1}+1_{\alpha\in(1,2)}.\label{Beta}
\end{align}
\end{enumerate}
\et

\section{$L^p$-maximal solution of linear nonlocal parabolic equation}

In the remainder of this paper, we shall fix $\alpha\in(0,2)$ and $m\in\mN$, and consider the following measurable maps: 
\begin{align*}
\mR_+\times\mR^m\ni(t,a)&\mapsto \nu_{t,a}\in\bL,\ \sigma_{t,a}\in\mM^d,\\
\mR_+\times\mR^d\ni(t,x)&\mapsto a_{t,x}\in\mR^m,\  b_{t,x}\in\mR^d.
\end{align*}
For a function $f:\mR_+\times\mR^d\to\mR^k$, where $k\in\mN$, the continuous modulus function associated to $f$ is defined by
$$
\hbar_f(\eps):=\sup_{|x-x'|\leq\eps}\sup_{t\geq 0}|f(t,x)-f(t,x')|.
$$
We make the following assumptions:
\begin{enumerate}[{\bf (H$_A$)}]
\item $\sigma_{t,a(t,\cdot)}$ and $ b$ are bounded, and for some $\nu^{(\alpha)}_1,\nu^{(\alpha)}_2,\nu^{(\alpha)}_3\in\bL^{(\alpha)}_{non}$,
$$
\nu^{(\alpha)}_1\leq \nu_{t,a}\leq\nu^{(\alpha)}_2, \ \ 1_{\alpha=1}\int_{r<|y|<R}y\nu_{t,a}(\dif y)=0,\ \ 0<r<R<\infty,
$$
and for some $\gamma_\sigma,\gamma_\nu\in(0,1)$ and $\kappa_0>0$,
\begin{align}
\hbar_\sigma(\eps)\leq \kappa_0\eps^{\gamma_\sigma}, \quad |\nu_{t,a}-\nu_{t,a'}|\leq |a-a'|^{\gamma_\nu}\nu^{(\alpha)}_3,\label{Ep12}
\end{align}
and $\lim_{\eps\to 0}(\hbar_a(\eps)+\hbar_b(\eps))=0$, and
\begin{align}\label{Sig}
\kappa_1:=\inf_{t\geq 0}\inf_{x\in\mR^d}\inf_{|\theta|=1}|\sigma_{t,a(t,x)}\theta|>0.
\end{align}
\end{enumerate}

Consider the following operator
$$
\sL u(x):=\sL_tu(x):=A_{t}u(x)+B_{t}u(x),
$$
where $A_t u(x):=A_{t,x}u(x)$ with
\begin{align}
A_{t,z}u(x):=\cL^{\nu_{t,a(t,z)}}_{\sigma_{t,a(t,z)}} u(x)+1_{\alpha=1}  b_{t,z}\cdot\nabla u(x),\label{AB1}
\end{align}
and $B_t$ is an abstract linear operator from $H^{\alpha,p}$ to $L^p$ and satisfies that
\begin{enumerate}[{\bf (H$^{ p_0}_B$)}]
\item For some $ p_0>1$ and any $p\geq  p_0$ and $\delta>0$, there exists a constant $C_\delta>0$ such that for all $u\in H^{\alpha,p}$ and $t\geq 0$,
\begin{align}
\|B_tu\|_p\leq \delta\|\Delta^{\frac{\alpha}{2}}u\|_p+C_\delta\|u\|_p.\label{Es10}
\end{align}
\end{enumerate}
Here $A_{t}$ is the principal part of $\sL$, and $B_{t}$ is a lower order perturbation term. 
The reason of introducing the extra function $a$ can be seen from the following lemma, 
and the following examples should be kept in mind.
\begin{align}\tag{\bf Ex}
\left\{
\begin{array}{ll}
\mbox{Let $m=d$, $a(t,x)=x$, $\sigma_{t,a(t,x)}=\sigma_{t,x}$ and $\nu_{t,a(t,x)}=\nu_{t,x}$.}\\
\mbox{Let $m=d^2$, $\nu_{t,a}=\nu^{(\alpha)}$ for some $\nu^{(\alpha)}\in\bL^{(\alpha)}_{non}$ and $\sigma_{t,a}=a$. }\\
\mbox{Let $m=1$, $\nu_{t,a}=a\nu^{\alpha}$ for some $\nu^{(\alpha)}\in\bL^{(\alpha)}_{non}$ and $\sigma_{t,a}=\sigma$. }
\end{array}
\right\}
\end{align}
\bl\label{Le31}
For $\eps>0$, let $\chi_\eps$ be a bounded measurable function with support in $D_\eps:=\big\{x\in\mR^d: |x|\leq\eps\big\}$.
Under {\bf (H$_A$)}, for any $p>\frac{m}{\gamma_\sigma(\alpha\wedge 1)\wedge\gamma_\nu}$ and 
$\beta\in(\frac{m}{p\gamma_\sigma},\alpha\wedge 1)$, 
there is a constant $C>0$ such that for all $u\in\mH^{\alpha,p}$, $\eps\in(0,1)$ and $(t,z)\in\mR_+\times\mR^d$,
$$
\big\|\big(\cL^{\nu_{t,a(t,\cdot)}}_{\sigma_{t,a(t,\cdot)}}u-\cL^{\nu_{t,a(t,z)}}_{\sigma_{t,a(t,z)}}u\big)\chi_\eps(\cdot-z)\big\|_p
\leq C\big(\hbar_a(\eps)\big)^{\gamma_\sigma\beta\wedge\gamma_\nu}\|\chi_\eps\|_\infty\|\Delta^{\frac{\alpha}{2}}u\|_p.
$$
\el
\begin{proof}
Fix $(t,z)\in\mR_+\times\mR^d$ and $\eps\in(0,1)$. Let $\beta\in(\frac{m}{p\gamma_\sigma},\alpha\wedge 1)$.
By {\bf (H$_A$)} and \eqref{EW88}, we have for all $a_1,a_2\in\mR^m$ with $|a_i-a(t,z)|\leq \hbar_a(\eps)$,
$$
\|\cL^{\nu_{t,a_1}}_{\sigma_{t,a_1}}u-\cL^{\nu_{t,a_2}}_{\sigma_{t,a_2}}u\|_p
\leq C_1|a_1-a_2|^{\gamma_\sigma\beta\wedge\gamma_\nu}\|\Delta^\frac{\alpha}{2}u\|_p,
$$
where $C_1$ is independent of $(t,z)$ and $\eps$.
Hence, by Lemma \ref{Le27} with $f(a, x):=\cL^{\nu_{t,a}}_{\sigma_{t,a}}u(x)$, we have
\begin{align*}
&\big\|\big(\cL^{\nu_{t,a(t,\cdot)}}_{\sigma_{t,a(t,\cdot)}}u-\cL^{\nu_{t,a(t,z)}}_{\sigma_{t,a(t,z)}}u\big)\chi_\eps(\cdot-z)\big\|_p\\
&\leq\Big\|\sup_{|y-z|\leq \eps}\big|\cL^{\nu_{t,a(t,y)}}_{\sigma_{t,a(t,y)}}u-\cL^{\nu_{t,a(t,z)}}_{\sigma_{t,a(t,z)}}u\big|\Big\|_p\|\chi_\eps\|_\infty\\
&\leq\Big\|\sup_{|a-a(t,z)|\leq\hbar_a(\eps)}\big|\cL^{\nu_{t,a}}_{\sigma_{t,a}}u-\cL^{\nu_{t,a(t,z)}}_{\sigma_{t,a(t,z)}}u\big|\Big\|_p\|\chi_\eps\|_\infty\\
&\leq C_2\big(\hbar_a(\eps)\big)^{\gamma_\sigma\beta\wedge\gamma_\nu}\|\Delta^{\frac{\alpha}{2}}u\|_p\|\chi_\eps\|_\infty.
\end{align*}
The proof is complete.
\end{proof}

For the simplicity of notation, we shall write
$$
\mL^p(T):=\mL^p(0,T),\ \ \mH^{\alpha,p}(T):=\mH^{\alpha,p}(0,T)
$$
and
$$
\mU^{\alpha, p}(T):=\mH^{\alpha,p}(T)\cap\big \{\p_tu\in\mL^p(T)\big\}.
$$
The aim of this section is to prove that
\bt\label{Main1}
Suppose that {\bf (H$_A$)} and {\bf (H$^{ p_0}_B$)} hold.  Let  $p>\frac{d}{\alpha\wedge 1}\vee\frac{m}{\gamma_\sigma(\alpha\wedge 1)\wedge\gamma_\nu}\vee p_0$ and $T>0$.
For any $\lambda\geq 0$ and $f\in \mL^p(T)$, there exists a unique $u_\lambda\in \mU^{\alpha,p}(T)$ such that for all $t\in[0,T]$,
\begin{align}
u_\lambda(t)=\int^t_0(\sL_s-\lambda) u_\lambda(s)\dif s+\int^t_0f(s)\dif s\ \mbox{ in $L^p$}. \label{EE5}
\end{align}
Moreover, we have
\begin{align}
\|u_\lambda(t)\|^p_p\leq C_p\Big(\tfrac{1}{\lambda}\wedge t\Big)^{p-1}\int^t_0\|f(s)\|^p_p\dif s,\label{Es8}
\end{align}
and
\begin{align}
\|u_\lambda\|_{\mU^{\alpha,p}(T)}\leq C_p\|f\|_{\mL^p(T)},\label{Es88}
\end{align}
where the constant $C_p$ is independent of $\lambda$.
\et
\subsection{Case of $B=0$ and $\nu_{t,a}=\nu_t$, $\sigma_{t,a}=\sigma_t$, $ b_{t,x}= b_t$}
In this subsection, we first consider the case of constant coefficients.
Let $N(\dif t,\dif y)$ be the Poisson random measure with intensity measure $\nu_{t}(\dif y)\dif t$. Let
$\tilde N(\dif t,\dif y):=N(\dif t,\dif y)-\nu_t(\dif y)\dif t$ be the compensated random
martingale measure. For $t\geq 0$, define
\begin{align}\label{XX}
\begin{split}
X_t:=1_{\alpha=1}&\int^t_0 b_r\dif r+\int^t_0\!\!\!\int_{\mR^d}\sigma_r y^{(\alpha)}\tilde N(\dif r,\dif y)\\
&+1_{\alpha\in(0,1)}\int^t_0\!\!\!\int_{\mR^d}\sigma_r yN(\dif r,\dif y).
\end{split}
\end{align}
For $\varphi\in C^2_b(\mR^d)$, by It\^o's formula we have
\begin{align*}
\mE\varphi(x+X_t-X_s)&=\varphi(x)+1_{\alpha=1}\mE\int^t_s b_r\cdot\nabla\varphi(x+X_r-X_s)\dif r\\
&+\mE\int^t_s\!\!\!\int_{\mR^d}\cJ^{(\alpha)}_\varphi(x, X_r-X_s+\sigma_r y)\nu_r(\dif y)\dif r.
\end{align*}
Thus, if we let
\begin{align}
\cT_{t,s}\varphi(x):=\mE\varphi\left(x+X_t-X_s\right),\label{EY2}
\end{align}
then
$$
\cT_{t,s}\varphi(x)=\varphi(x)+\int^t_sA_r\cT_{r,s}\varphi(x)\dif r.
$$

The following result is a simple application of \cite[Theorem 4.2]{Zh1}.
\bt\label{Th1}
For and $T>0$, $p>1$ and $f\in \mL^p(T)$, let
$$
u_\lambda(t,x):=\int^t_0 \mathrm{e}^{-\lambda (t-s)}\cT_{t,s}f(s,x)\dif s.
$$
Under {\bf (H$_A$)},  $u_\lambda$ is the unique solution of equation (\ref{EE5})  with
\begin{align}\label{ER10}
\|u_\lambda(t)\|_p^p\leq \left(\frac{1-\mathrm{e}^{-\lambda t}}{\lambda}\right)^{p-1}\int^t_0\mathrm{e}^{-\lambda(t-s)}\|f(s)\|^p_p\dif s,
\end{align}
and for some $C=C\big(p,d,\alpha,\kappa_0,\nu_1^{(\alpha)}, \nu_2^{(\alpha)}\big)>0$,
\begin{align}
\|u_\lambda\|_{\mU^{\alpha,p}(T)}\leq C\|f\|_{\mL^p(T)}.\label{ER11}
\end{align}
\et
\begin{proof}
It suffices to prove estimates (\ref{Es8}) and (\ref{Es88}). By a mollifying technique, we may assume that $f\in L^p([0,T];\cap_{\beta\geq 0} H^{\beta,p})$. 

Let $N^{(1)}(\dif t,\dif y)$ and $N^{(2)}(\dif t,\dif y)$ be two independent Poisson random measures with
intensity measures $\nu^{(\alpha)}_1(\dif y)\dif t$ and $(\nu_t(\dif y)-\nu^{(\alpha)}_1(\dif y))\dif t$ respectively, where $\nu^{(\alpha)}_1$ 
is the lower bound of $\nu_t$ from {\bf (H$_A$)}.
Let $X^{(2)}_t$ be defined by (\ref{XX}) in terms of $N^{(2)}$, and $X^{(1)}_t$ be defined by
$$
X^{(1)}_t:=\int^t_0\!\!\!\int_{\mR^d}\sigma_r y^{(\alpha)}\tilde N^{(1)}(\dif r,\dif y)
+1_{\alpha\in(0,1)}\int^t_0\!\!\!\int_{\mR^d}\sigma_r yN^{(1)}(\dif r,\dif y).
$$
Set for $\varphi\in C^2_b(\mR^d)$,
$$
\cT^{(i)}_{t,s}\varphi(x):=\mE \varphi\left(x+X^{(i)}_t-X^{(i)}_s\right),\ \ i=1,2.
$$
Since $X^{(1)}_\cdot$ and $X^{(2)}_\cdot$ are independent and
$$
X^{(1)}_\cdot+X^{(2)}_\cdot\stackrel{(d)}{=}X_\cdot,
$$
we have
\begin{align*}
\cT_{t,s}\varphi(x)=\cT^{(2)}_{t,s}\cT^{(1)}_{t,s}\varphi(x)
=\mE\cT^{(1)}_{t,s}\varphi(x+X^{(2)}_t-X^{(2)}_s).
\end{align*}
Thus, by Jensen's inequality and \cite[Theorem 4.2]{Zh1}, there exits a constant $C=C(d,p,\alpha,\kappa_0,\nu^{(\alpha)}_1)>0$ such that
\begin{align}
\int^T_0\|\Delta^{\frac{\alpha}{2}} u_\lambda(t)\|_p^p\dif t
&\leq\mE\int^T_0\left\|\Delta^{\frac{\alpha}{2}} \int^t_0\mathrm{e}^{-\lambda (t-s)}\cT^{(1)}_{t,s}f(s,\cdot+X^{(2)}_t-X^{(2)}_s)\dif s\right\|_p^p\dif t\no\\
&=\mE\int^T_0\left\|\Delta^{\frac{\alpha}{2}}\int^t_0\mathrm{e}^{-\lambda (t-s)}\cT^{(1)}_{t,s}f(s,\cdot-X^{(1)}_s)\dif s\right\|_p^p\dif t\no\\
&\leq C\mE\int^T_0\left\|f(s,\cdot-X^{(1)}_s)\right\|_p^p\dif s=C\int^T_0\|f(s)\|^p_p\dif s.\label{ER09}
\end{align}
On the other hand, by Minkowskii and H\"older's inequalities,  we have
\begin{align*}
\|u(t)\|_p^p&\leq \left(\int^t_0\mathrm{e}^{-\lambda(t-s)}\|\cT_{t,s}f(s)\|_p\dif s\right)^p\leq \left(\int^t_0\mathrm{e}^{-\lambda(t-s)}\|f(s)\|_p\dif s\right)^p\\
&\leq \left(\frac{1-\mathrm{e}^{-\lambda t}}{\lambda}\right)^{p-1}\int^t_0\mathrm{e}^{-\lambda(t-s)}\|f(s)\|^p_p\dif s,
\end{align*}
which then gives (\ref{ER10}). Moreover, by Fubini's theorem, we also have
\begin{align*}
\int^T_0\|u(t)\|_p^p\dif t\leq\frac{1}{\lambda^p}\int^T_0\|f(s)\|^p_p\dif s,
\end{align*}
which together with (\ref{EE5}), (\ref{Es17}) and (\ref{ER09}) yields 
$$
\|\p_t u\|_{\mL^p(T)}\leq\|A u\|_{\mL^p(T)}+\lambda\|u\|_{\mL^p(T)}+\|f\|_{\mL^p(T)}\leq C\|f\|_{\mL^p(T)}.
$$
The proof is complete.
\end{proof}

\subsection{Freezing function and auxiliary estimates}
Let $p\geq 1$ and $\phi\in C^\infty_c(\mR^d)$ be a nonnegative symmetric function with support in the unit ball and satisfy
$$
\int_{\mR^d}\phi^p(x)\dif x=1.
$$
For $\delta\in(0,1)$, let us set
\begin{align}
\phi^z_\delta(x):=\delta^{-d/p}\phi(\delta^{-1}(x-z)).\label{zeta}
\end{align}
Then 
\begin{align}
\|\phi_\delta^\cdot (x)\|_p^p=1,\ \ \forall x\in\mR^d,\label{zeta1}
\end{align}
and $\{\phi_\delta^z(\cdot),\delta\in(0,1),z\in\mR^d\}$ will serve as a family of freezing functions as shown in the following two crucial lemmas.

\bl\label{Le9}
For $p\geq 1$ and $\delta\in(0,1)$, there exist two constants $C_1, C_2>0$ depending only on $d, p,\alpha,\delta$ and $\phi$ 
such that for all $u\in H^{\alpha,p}$,
\begin{align}\label{Es5}
\tfrac{1}{2}\|\Delta^{\frac{\alpha}{2}}u\|_p-C_1\|u\|_p\leq\left(\int_{\mR^d}\|\Delta^{\frac{\alpha}{2}}(u\phi_\delta^z)\|_p^p\dif z\right)^{\frac{1}{p}}
\leq \tfrac{3}{2}\|\Delta^{\frac{\alpha}{2}}u\|_p+C_2\|u\|_p.
\end{align}
\el
\begin{proof}
By (\ref{Es7}) and (\ref{For}), we have
\begin{align}\label{Es3}
\begin{split}
I_\delta^z(x)&:=\Delta^{\frac{\alpha}{2}}(u\phi_\delta^z)(x)-\phi_\delta^z(x)\Delta^{\frac{\alpha}{2}}u(x)-u(x)\Delta^{\frac{\alpha}{2}}\phi_\delta^z(x)\\
&=c^{-1}_{d,\alpha}\int_{\mR^d}(u(x+y)-u(x))(\phi_\delta^z(x+y)-\phi_\delta^z(x))\frac{\dif y}{|y|^{d+\alpha}}.
\end{split}
\end{align}
By definitions, it is easy to see that
\begin{align}
\int_{\mR^d}|\Delta^{\frac{\alpha}{2}}\phi_\delta^z(x)|^p\dif z=\int_{\mR^d}|\Delta^{\frac{\alpha}{2}}\phi_\delta(z)|^p\dif z<\infty,\label{ER1}
\end{align}
and
\begin{align}
\sup_{x\in\mR^d}\left(\int_{\mR^d}|\phi_\delta^z(x+y)-\phi_\delta^z(x)|^p\dif z\right)^{\frac{1}{p}}\leq (C_\delta|y|)\wedge 2.\label{ER2}
\end{align}
Moreover, for any $\beta\in(0,1)$, by (\ref{Es34}) we also have
\begin{align}
\|u(\cdot+y)-u(\cdot)\|_p\leq C_{d,\beta}|y|^\beta\|\Delta^{\frac{\beta}{2}}u\|_p.\label{ER3}
\end{align}
Hence, for any $\beta\in(0,1\wedge \alpha)$, using (\ref{ER1}), (\ref{ER2}) and (\ref{ER3}), and by Minkowskii's inequality and interpolation inequality (\ref{Inter}), we derive that
\begin{align}\label{ER4}
\begin{split}
\left(\int_{\mR^d}\!\!\int_{\mR^d}|I_\delta^z(x)|^p\dif x\dif z\right)^{\frac{1}{p}}
&\leq C_\delta\|u\|_p+C_\delta \|\Delta^{\frac{\beta}{2}}u\|_p
\leq \tfrac{1}{2}\|\Delta^{\frac{\alpha}{2}}u\|_p+C_{\delta}\|u\|_p.
\end{split}
\end{align}
Substituting this into (\ref{Es3}) and using (\ref{zeta1}), we obtain (\ref{Es5}).
\end{proof}

\bl\label{Le34}
Under {\bf (H$_A$)}, for any $p>\frac{d}{\alpha\wedge 1}\vee\frac{m}{\gamma_\sigma(\alpha\wedge 1)\wedge \gamma_\nu}$,
there is a function $\ell(\delta)$ with $\lim_{\delta\to 0}\ell(\delta)=0$ such that for all $u\in H^{\alpha,p}$ and $t\geq 0$,
\begin{align}\label{ER5}
\begin{split}
\left(\int_{\mR^d}\|(A_{t}u)\phi_\delta^z-A_{t,z}(u\phi_\delta^z)\|^p_p\dif z\right)^{1/p}\leq \ell(\delta)\|\Delta^{\frac{\alpha}{2}}u\|_p+C_\delta\|u\|_p.
\end{split}
\end{align}
In particular,
\begin{align}
\|A_tu\|_p\leq C\|u\|_{ \alpha,p}.\label{Es17}
\end{align}
\el
\begin{proof}
In the following we shall drop the time variable since it does not play any role in the proof.
First of all, by \eqref{zeta1}, it is easy to see that
\begin{align}
&\int_{\mR^d}\|( b_{\cdot}\cdot\nabla u)\phi_\delta^z- b_{z}\cdot\nabla (u\phi_\delta^z)\|_p^p\dif z
\preceq\int_{\mR^d}\|(( b_{\cdot}- b_z)\cdot\nabla u)\phi_\delta^z\|^p_p\dif z\no\\
&\quad+\| b\|^p_\infty\int_{\mR^d}\|u\nabla\phi_\delta^z\|_p^p\dif z
\leq\hbar_b(\delta)^p\|\nabla u\|^p_p+\delta^{-p}\| b\|^p_\infty\|u\|^p_p\|\nabla\phi\|^p_p.\label{LK1}
\end{align}

Below, for the simplicity of notation, we write
$$
\mu_x:=\nu_{a(x)}, \ \ \Theta_x:=\sigma_{a(x)}.
$$
Let $\chi:\mR^d\to[0,1]$ be a smooth function with $\chi(x)=1$ for $|x|<2$ and $\chi(x)=0$ for $|x|>4$. 
For $\delta\in(0,1)$, let $\chi^z_\delta(x):=\chi(\delta^{-1}(x-z))$.
Let us write 
\begin{align*}
\phi_\delta^z(x)\cL^{\mu_x}_{\Theta_x} u(x)-\cL^{\mu_z}_{\Theta_z}(u\phi_\delta^z)(x)=I^z_1(x)+I^z_2(x)-I^z_3(x),
\end{align*}
where
\begin{align*}
I^z_1(x)&:=\big(\cL^{\mu_x}_{\Theta_x}(u\phi_\delta^z)(x)-\cL^{\mu_z}_{\Theta_z}(u\phi_\delta^z)(x)\big)\chi_\delta^z(x),\\
I^z_2(x)&:=\big(\cL^{\mu_x}_{\Theta_x}(u\phi_\delta^z)(x)-\cL^{\mu_z}_{\Theta_z}(u\phi_\delta^z)(x)\big))(1-\chi_\delta^z(x)),\\
I^z_3(x)&:=u(x)\cL^{\mu_x}_{\Theta_x}\phi_\delta^z(x)+\int_{\mR^d}\big(u(x+\Theta_xy)-u(x)\big)\\
&\qquad\qquad\times\big(\phi_\delta^z(x+\Theta_x y)-\phi_\delta^z(x)\big)\mu_x(\dif y).
\end{align*}
For $I_1^z(x)$, by Lemma \ref{Le31} and \eqref{Es5}, we have for $\beta\in(0,\alpha\wedge 1)$,
\begin{align}\label{LK3}
\begin{split}
\left(\int_{\mR^d}\|I^z_1\|_p^p\dif z\right)^{1/p}
&\leq C\big(\hbar_a(4\delta)\big)^{\gamma_\sigma\beta\wedge\gamma_\nu}\left(\int_{\mR^d}\|\Delta^{\frac{\alpha}{2}}(u\phi^z_\delta)\|_p^p\dif z\right)^{1/p}\\
&\leq C\big(\hbar_a(4\delta)\big)^{\gamma_\sigma\beta\wedge\gamma_\nu}\|\Delta^{\frac{\alpha}{2}}u\|_p+C_\delta\|u\|_p.
\end{split}
\end{align}
Noticing that
$$
\chi_\delta^z(x)=1,\ \ |x-z|\leq 2\delta,\ \ \phi_\delta^z(x)=0,\ \ |x-z|\geq\delta,
$$
by the definition of $\cL^{\mu_x}_{\Theta_x}$ and $1_{\alpha=1}\int_{r<|y|<R}y\mu_{x}(\dif y)=0$, it is easy to see that
\begin{align*}
\|I^z_2\|^p_p&\leq\int_{\mR^d}\left|\int_{|y|>\frac{\delta}{\kappa}}
\cJ^{(\alpha)}_{u\phi_\delta^z}(x,\Theta_xy)\mu_x(\dif y)-\int_{|y|>\frac{\delta}{\kappa}}\cJ^{(\alpha)}_{u\phi_\delta^z}(x,\Theta_zy)\mu_z(\dif y)\right|^p\dif x\\
&=\int_{\mR^d}\left|\int_{|y|>\frac{\delta}{\kappa}}
\cJ^{(\theta)}_{u\phi_\delta^z}(x,\Theta_xy)\mu_x(\dif y)-\int_{|y|>\frac{\delta}{\kappa}}\cJ^{(\theta)}_{u\phi_\delta^z}(x,\Theta_zy)\mu_z(\dif y)\right|^p\dif x,
\end{align*}
where $\kappa:=\|\sigma_a\|_\infty$, and $\theta\in(\frac{d}{p},\alpha)$ for $\alpha\in(0,1]$ and $\theta\in(1,\alpha)$ for $\alpha\in(1,2)$. Noticing that
\begin{align*}
\int_{\mR^d}\left|\sup_{y\not=0}\Big(|y|^{-\theta}|\cJ^{(\theta)}_{u\phi_\delta^z}(x,\Theta_x y)|\Big)\right|^p\dif x
\leq \|\sigma_a\|_\infty^{p\theta}\int_{\mR^d}\left|\sup_{y\not=0}\Big(|y|^{-\theta}|\cJ^{(\theta)}_{u\phi_\delta^z}(x,y)|\Big)\right|^p\dif x,
\end{align*}
by (\ref{Es34}), \eqref{Es5} and (\ref{Inter}), we have
\begin{align}
\left(\int_{\mR^d}\|I^z_2\|^p_p\dif z\right)^{1/p}&
\preceq\left(\int_{\mR^d}\|\Delta^{\frac{\theta}{2}}(u\phi_\delta^z)\|_p^p\dif z\right)^{1/p}\int_{|y|>\frac{\delta}{\kappa}}|y|^\theta\nu^{(\alpha)}_2(\dif y)\no\\
&\preceq C_\delta\big(\|\Delta^{\frac{\theta}{2}}u\|_p+\|u\|_p\big)\leq\delta\|\Delta^{\frac{\alpha}{2}}u\|_p+C_\delta\|u\|_p.\label{LK2}
\end{align}
On the other hand, by definition, it is easy to see that
$$
\sup_{x}\int_{\mR^d}|\cL^{\mu_x}_{\Theta_x}\phi_\delta^z(x)|^p\dif z\leq C_\delta
$$
and
$$
\sup_{x\in\mR^d}\left(\int_{\mR^d}|\phi_\delta^z(x+\Theta_x y)-\phi_\delta^z(x)|^p\dif z\right)^{\frac{1}{p}}\leq (C_\delta|y|)\wedge 2.
$$
Moreover, for any $\beta\in(\frac{d}{p},\alpha\wedge 1)$, by (\ref{Es34}) we also have
$$
\|u(\cdot+\Theta_\cdot y)-u(\cdot)\|_p\leq C_{d,\beta}|y|^\beta\|\Delta^{\frac{\beta}{2}}u\|_p.
$$
For $I^\delta_3(z,x)$, as in  estimating (\ref{ER4}), we have
\begin{align}\label{LK4}
\left(\int_{\mR^d}\|I^z_3\|^p_p\dif z\right)^{1/p}\leq \delta\|\Delta^{\frac{\alpha}{2}}u\|_p+C_\delta\|u\|_p.
\end{align}
Combining \eqref{LK1}-\eqref{LK4}, we obtain (\ref{ER5}).

Finally, by \eqref{zeta1}, \eqref{EW888}, \eqref{ER5} and \eqref{Es5}, we have
\begin{align*}
\|A_tu\|_p&=\left(\int_{\mR^d}\|(A_t u)\phi^z_\delta\|_p^p\dif z\right)^{1/p}\leq
\left(\int_{\mR^d}\|A_{t,z}(u\phi^z_\delta)\|_p^p\dif z\right)^{1/p}\\
&\qquad+\left(\int_{\mR^d}\|(A_t u)\phi^z_\delta-A_{t,z}(u\phi^z_\delta)\|_p^p\dif z\right)^{1/p}\\
&\preceq \left(\int_{\mR^d}\|\Delta^{\frac{\alpha}{2}}(u\phi^z_\delta)\|_p^p\dif z\right)^{1/p}+\|u\|_{\alpha,p}\preceq\|u\|_{\alpha,p}.
\end{align*}
The proof is complete.
\end{proof}

\subsection{Proof of Theorem \ref{Main1}}
We divide the proof into two steps.
\\
\\
(1)  We first prove the a priori estimates (\ref{Es8}) and (\ref{Es88}). 
Let $\varrho:\mR^d\to\mR_+$ be a smooth function with support in the unit ball and $\int\varrho=1$. For $\eps\in(0,1)$, let
$\varrho_\eps(x):=\eps^{-d}\varrho(\eps^{-1}x).$
Taking convolutions for both sides of (\ref{EE5}) with respect to $\varrho_\eps$, we have
\begin{align}
\p_tu^\eps_\lambda=(\sL-\lambda) u^\eps_\lambda+h_\eps,\label{Eq333}
\end{align}
where $u^\eps_\lambda:=u_\lambda*\varrho_\eps$ and
$$
h_\eps:=f*\varrho_\eps+(\sL u_\lambda)*\varrho_\eps-\sL (u_\lambda*\varrho_\eps).
$$
By (\ref{Es17}), {\bf (H$^{p_0}_B$)} and the property of convolutions, we have
$$
\lim_{\eps\to0}\int^T_0\|h_\eps(t)-f(t)\|^p_p\dif t=0.
$$
Below, we use the method of freezing the coefficients to prove that for all $T\geq 0$ and $t\in[0,T]$,
\begin{align}
\|u^\eps_\lambda(t)\|^p_p\leq C\Big(\tfrac{1}{\lambda}\wedge t\Big)^{p-1}\int^t_0\|h_\eps(s)\|_p^p\dif s,\ \ \
\|u^\eps_\lambda\|_{\mU^{\alpha,p}(T)}\leq C\|h_\eps\|_{\mL^p(T)},\label{Ep7}
\end{align}
where the constant $C$ is independent of $\eps$ and $\lambda$. After proving this estimate, (\ref{Es8}) and (\ref{Es88})
immediately follows by Fatou's lemma and taking limits for (\ref{Ep7}).

Let $\phi^\delta_z$ be defined by (\ref{zeta}). For the simplicity of notation, we drop the subscript $\eps,\lambda$ and $\delta$ below. Multiplying both sides of (\ref{Eq333}) by $\phi_z$, we have
$$
\p_t(u\phi_z)=(A_{t,z}-\lambda) (u\phi_z)+g^\phi_z,
$$
where
$$
g^\phi_z:=(Au)\phi_z-A_{t,z} (u\phi_z)+(Bu+h)\phi_z.
$$
By Lemma \ref{Le9} and Theorem \ref{Th1}, we have
\begin{align}
\|\Delta^{\frac{\alpha}{2}}u\|_{\mL^p(T)}&\leq\tfrac{3}{2}\left(\int_{\mR^d}\|\Delta^{\frac{\alpha}{2}}(u\phi_z)\|^p_{\mL^p(T)}\dif z\right)^{1/p}
+\|u\|_{\mL^p(T)}\no\\
&\leq C\left(\int_{\mR^d}\|g^\phi_z\|^p_{\mL^p(T)}\dif z\right)^{1/p}+C_\delta\|u\|_{\mL^p(T)}.\label{Es6}
\end{align}
Recalling $\|\phi_\cdot(x)\|^p_p=1$, by definitions, Lemmas \ref{Le34} and {\bf (H$^{p_0}_B$)}, we have
\begin{align}
\left(\int_{\mR^d}\|g^\phi_z(t)\|^p_p\dif z\right)^{1/p}&\leq
\left(\int_{\mR^d}\|(A_tu)\phi_z-A_{t,z} (u\phi_z)\|^p_p\dif z\right)^{1/p}+\|(Bu+h)(t)\|_p\no\\
&\leq (\ell(\delta)+\delta)\|\Delta^{\frac{\alpha}{2}} u(t)\|_p+C_{\delta}\|u(t)\|_p+\|h(t)\|_p.\label{Es11}
\end{align}
Substituting this into (\ref{Es6}) and letting $\delta$ be small enough, we obtain
\begin{align}
\|\Delta^{\frac{\alpha}{2}}u\|_{\mL^p(T)}\leq C\left(\|u\|_{\mL^p(T)}+\|h\|_{\mL^p(T)}\right).\label{Es12}
\end{align}
On the other hand, by Theorem \ref{Th1}, we also have for all $t\in[0,1]$,
$$
\|u(t)\|^p_p=\int_{\mR^d}\|(u\phi_z)(t)\|^p_p\dif z\leq C\Big(\tfrac{1}{\lambda}\wedge t\Big)^{p-1}\int^t_0\!\!\! \int_{\mR^d}\|g^\phi_z(s)\|^p_p\dif z\dif s,
$$
which together with (\ref{Es11}), (\ref{Es12}) gives
$$
\|u(t)\|^p_p\leq C\Big(\tfrac{1}{\lambda}\wedge t\Big)^{p-1}\left(\int^t_0(\|u(s)\|^p_p+\|h(s)\|^p_p)\dif s\right).
$$
By Gronwall's inequality, we obtain
\begin{align}
\|u(t)\|^p_p\leq C\Big(\tfrac{1}{\lambda}\wedge t\Big)^{p-1}\int^t_0 \|h(s)\|^p_p\dif s,\label{EI1}
\end{align}
which together with (\ref{Es12}) yields (\ref{Ep7}).
\\
\\
(2) In this step we use the classical continuity method to prove the existence of solutions (cf. \cite{Kr}). For $\tau\in[0,1]$, define an operator
$$
U_\tau:=\p_t-\tau(\sL-\lambda)-(1-\tau)\cL^{\nu^{(\alpha)}_1}.
$$
By (\ref{Es17}) and (\ref{Es10}), it is easy to see that
\begin{align}
U_\tau: \mU^{\alpha,p}(T)\to \mL^p(T).\label{Eq6}
\end{align}
For $\tau=0$ and $f\in \mL^p(T)$, by Theorem \ref{Th1}, there is a unique $u\in\mU^{\alpha,p}(T)$ such that
$$
U_0 u=\p_t u-\cL^{\nu^{(\alpha)}_1} u=f.
$$
Suppose now that for some $\tau_0\in[0,1)$, and for any $f\in \mL^p(T)$, the equation
$$
U_{\tau_0}u=f
$$
admits a unique solution $u\in{\mU^{\alpha,p}(T)}$. Under this assumption, for fixed $f\in\mL^p(T)$ and $\tau\in[\tau_0,1]$, 
and for any $u\in{\mU^{\alpha,p}(T)}$, by (\ref{Eq6}), the equation
\begin{align}
U_{\tau_0}w=f+(U_{\tau_0}-U_\tau)u\label{Eq0}
\end{align}
admits a unique solution $w\in{\mU^{\alpha,p}(T)}$. Introduce an operator
$$
Q^f_\tau: u\mapsto w=Q^f_\tau u.
$$
We now use the apriori estimate (\ref{Es88}) to show that there exists an $\eps>0$
independent of $\tau_0$ such that for all $\tau\in[\tau_0,\tau_0+\eps]$,
$$
Q^f_\tau:{\mU^{\alpha,p}(T)}\to{\mU^{\alpha,p}(T)}
$$
is a contraction operator.
Let $u_1,u_2\in{\mU^{\alpha,p}(T)}$ and $w_i=Q^f_{\tau}u_i,i=1,2$. By equation (\ref{Eq0}), we have
\begin{align*}
U_{\tau_0}(w_1-w_2)=(U_{\tau_0}-U_\tau)(u_1-u_2)=(\tau_0-\tau)((\sL-\lambda)-\cL^{\nu^{(\alpha)}_1})(u_1-u_2).
\end{align*}
By (\ref{Es88}), (\ref{Es17}) and (\ref{Es10}), it is not hard to see that
\begin{align*}
&\|Q^f_\tau u_1-Q^f_\tau u_2\|_{{\mU^{\alpha,p}(T)}}=\|w_1-w_2\|_{\mU^{\alpha,p}(T)}\\
&\leq C|\tau_0-\tau|\cdot\|((\sL-\lambda)-\cL^{\nu^{(\alpha)}_1})(u_1-u_2)\|_{\mL^p(T)}\\
&\leq C_0|\tau_0-\tau|\cdot\|u_1-u_2\|_{{\mU^{\alpha,p}(T)}},
\end{align*}
where $C_0$ is independent of $\tau,\tau_0$ and $u_1,u_2, f$. Taking $\eps=1/(2C_0)$, one sees that for all $\tau\in[\tau_0,\tau_0+\eps]$,
$$
Q^f_\tau: {\mU^{\alpha,p}(T)}\to {\mU^{\alpha,p}(T)}
$$
is a $1/2$-contraction operator. By the fixed point theorem, for each $\tau\in[\tau_0,\tau_0+\eps]$, there exists a unique $u\in{\mU^{\alpha,p}(T)}$ such that
$$
Q^f_\tau u=u,
$$
which means that
$$
U_\tau u=f.
$$
Now starting from $\tau=0$, after repeating the above construction
$[\frac{1}{\eps}]+1$-steps, one obtains that for any $f\in\mL^p(T)$,
$$
U_1u=f
$$
admits a unique solution $u\in{\mU^{\alpha,p}(T)}$.

\section{Uniqueness of martingale solutions}

Let $\sM(\mR^d)$ be the set of all signed measures over $\mR^d$ endowed with weak convergence topology. 
In this section, we shall take $B$ as the following concrete form: for some $\bar\alpha\in(0,\alpha)$,
\begin{align}\label{AB2}
B_tu(x):=\int_{\mR^d}\cJ^{(\bar\alpha)}_u(x,\bar\sigma_{t,x}y)\bar\nu_{t,x}(\dif y)+1_{\alpha\in(1,2)}\bar b_{t,x}\cdot\nabla u,
\end{align}
where 
$$
\mR_+\times\mR^d\ni(t,x)\mapsto\bar\sigma_{t,x}\in\mM^d,\ \bar\nu_{t,x}\in\sM(\mR^d),\ \ \bar b_{t,x}\in\mR^d
$$
are Borel measurable and satisfy that
\begin{enumerate}[{\bf (H$'_B$)}]
\item $\bar\sigma$ and $\bar b$ are bounded, and $|\bar\nu_{t,x}|\leq \nu^{(\bar\alpha)}_4$ for some $\nu^{(\bar\alpha)}_4\in \bL^{(\bar\alpha)}_{non}$.
\end{enumerate}

The following lemma is direct by definition and Lemma \ref{Le22}.
\bl\label{Le41}
Under {\bf (H$'_B$)}, we have
\begin{align*}
&\left|\int_{\mR^d}\cJ^{(\bar\alpha)}_u(x,\bar\sigma_{t,x}y)\bar\nu_{t,x}(\dif y)\right|
\leq \|\bar\sigma\|^{\theta_1}_\infty\sup_{y\not=0}\frac{\cJ^{(\theta_1)}_u(x,y)}{|y|^{\theta_1}}\int_{|y|\leq 1}|y|^{\theta_1}\nu^{(\bar\alpha)}_4(\dif y)\\
&+\|\bar\sigma\|^{\theta_2}_\infty\sup_{y\not=0}\frac{\cJ^{(\theta_2)}_u(x,y)}{|y|^{\theta_2}}\int_{|y|>1}|y|^{\theta_2}\nu^{(\bar\alpha)}_4(\dif y)
+1_{\bar\alpha=1}\nu^{(\bar\alpha)}_4(|y|>1)|\nabla u(x)|,
\end{align*}
where $\theta_2<\bar\alpha<\theta_1$ are chosen in the following way:
$$
\left\{
\begin{aligned}
&\theta_1\in(\bar\alpha,\alpha\wedge 1),\ \ \theta_2\in(0,\bar\alpha),\quad \bar\alpha\in(0,1);\\
&\theta_1\in(1,\alpha),\ \ \theta_2\in(0,1),\quad \bar\alpha=1;\\
&\theta_1\in(\bar\alpha,\alpha),\ \ \theta_2\in(1,\bar\alpha),\quad \bar\alpha\in(1,2).
\end{aligned}
\right.
$$
In particular, {\bf (H$^{p_0}_B$)} holds with $p_0=\frac{d}{\bar\alpha}\vee 1$.
\el

Let $\Omega=D(\mR_+;\mR^d)$ be the space of all right continuous functions with left hand limits, which is endowed with the Skorokhod metric. Let 
$$
X_t(\omega):=\omega_t
$$
be the coordinate process on $\Omega$, and
$$
\sF_t:=\sigma\big\{X_s: s\in[0,t]\big\}, \ \ \sF:=\sF_\infty.
$$
\bd
\begin{enumerate}[(i)]
\item (Martingale solution) For fixed $(s,x)\in\mR_+\times\mR^d$, we say that a probability measure $\mP$ on $(\Omega,\sF)$ is a solution to the martingale problem for $\sL$
starting from $(s,x)$ if $\mP(X_r=x, r\in[0,s])=1$ and for all $\varphi\in C^\infty_c(\mR^d)$,
$$
t\mapsto \varphi(X_t)-\int^t_s\sL_{r}\varphi(X_r)\dif r=:M^\varphi_t,\ \ t\geq s,
$$
is an $\sF_t$-martingale under $\mP$. The set of all martingale solutions to the martingale problem for $\sL$ 
with starting point $(s,x)$ is denoted by $\Gamma^\sL_{s,x}$. 
\item (Krylov's type estimate) Let $\mP\in\Gamma^\sL_{s,x}$. One says that Krylov's type estimate holds for $\mP$ if for some $p_1>1$ and any $p\geq p_1$ and $T>s$,
there exits a constant $C>0$ such that for all $s\leq t_1\leq t_2\leq T$ and $f\in\mL^p(T)$,
$$
\mE\left(\int^{t_2}_{t_1} f(s,X_s)\dif s\Big|\sF_{t_1}\right)\leq C\|f\|_{\mL^p(T)}.
$$
All the martingale solutions with the above property is denoted 
by $\tilde\Gamma^\sL_{s,x}$.
\end{enumerate}
\ed

\br\label{Re}
Under {\bf (H$_A$)} and {\bf (H$'_B$)},
by suitable approximation, for any $u\in C(\mR_+;  H^\infty)$ with $\p_t u\in L^1_{loc}(\mR_+; H^\infty)$,
$$
t\mapsto u(t,X_t)-\int^t_s(\p_r+\sL_r)u(r,X_r)\dif r
$$
is still an $\sF_t$-martingale after time $s$.
\er

We now show the following important Krylov's type estimate.
\bt\label{Main}
Assume that {\bf (H$_A$)} and {\bf (H$'_B$)} hold.
For any $T>0$ and $p>\tfrac{m}{\gamma_\sigma(\alpha\wedge 1)\wedge\gamma_\nu}\vee\frac{d}{\bar\alpha}\vee \Big(\tfrac{d}{\alpha}+1\Big)\vee\frac{d}{\alpha\wedge 1}$, 
there exists a constant $C>0$ only depending on the bounds and parameters appearing in {\bf (H$_A$)} and {\bf (H$'_B$)}
such that for any $(s,x)\in[0,T]\times\mR^d$ and $\mP\in\Gamma_{s,x}^\sL$, 
any $s\leq t_1\leq t_2\leq T$ and $f\in\mL^p(s,T)$, if for some $\gamma_0\in(0,1)$,
\begin{align}\label{NG6}
\sup_{\eps\in(0,1)}\Big((\hbar_a(\eps)+\hbar_b(\eps))\eps^{-\gamma_0}\Big)<+\infty,
\end{align}
then it holds that for any $\beta\in(\frac{d}{p},\alpha(1-\frac{1}{p}))$,
\begin{align}
\mE\left(\int^{t_2}_{t_1}f(r,X_r)\dif r\Big|\sF_{t_1}\right)\leq C(t_2-t_1)^{1-\frac{\beta}{\alpha}-\frac{1}{p}}\|f\|_{\mL^p(s,T)},\label{Eq24}
\end{align}
where the expectation $\mE$ is taken with respect to $\mP$. 
\et
\begin{proof}
Below, without loss of generality, we assume $s=0$.
Let 
\begin{align}\label{AS1}
p>\tfrac{m}{\gamma_\sigma(\alpha\wedge 1)\wedge\gamma_\nu}\vee\tfrac{d}{\bar\alpha}\vee\Big(\tfrac{d}{\alpha}+1\Big)\vee\tfrac{d}{\alpha\wedge 1},
\end{align}
and
\begin{align}\label{AS2}
\theta:=1_{\alpha\in(0,1]}\tfrac{\alpha+\bar\alpha}{2}+1_{\alpha\in(1,2)}\tfrac{\alpha+\bar\alpha\vee 1}{2},\ \ 
q>\tfrac{d}{\gamma_0(\gamma_\sigma(\alpha\wedge 1)\wedge\gamma_\nu)}\vee \tfrac{d}{\theta-\bar\alpha}\vee\tfrac{\alpha}{\alpha-\theta}\vee p. 
\end{align}
We divide the proof into two steps.
\\
\\
(1) First of all, we prove the following estimate of  Krylov's type: for any $f\in \mL^q(T)$,
\begin{align}
\mE \left(\int^T_0 f(r,X_r)\dif r\right)\leq C\|f\|_{\mL^q(T)}.\label{EE89}
\end{align}
By a standard approximation, we may assume $f\in C_c((0,T)\times\mR^d)$.
By Theorem \ref{Main1}, there exists a unique solution $u\in {\mU^{\alpha,q}(T)}$ to 
\begin{align}
\p_t u+Au=f,\ \  u(T)=0,\label{EP1}
\end{align}
with
\begin{align}
\|u\|_{\mU^{\alpha,q}(T)}\leq C\|f\|_{\mL^q(T)}.\label{EP2}
\end{align}
Let
$$
u_\eps(t,x):=u(t)*\varrho_\eps(x),\ f_\eps(t,x):=f(t)*\varrho_\eps(x).
$$
By Remark \ref{Re} and equation \eqref{EP1}, we have
\begin{align}\label{EP5}
\begin{split}
-u_\eps(0,x)&=\mE\left(\int^T_0(\p_r u_\eps(r,X_r)+\sL_r u_\eps(r,X_r))\dif r\right)\\
&=\mE\left(\int^T_0f_\eps(r,X_r)\dif r\right)+\int^T_0\Big(\Lambda_\eps(r)+\Lambda'_\eps(r)\Big)\dif r,
\end{split}
\end{align}
where
\begin{align*}
\Lambda_\eps(r)&:=\mE\left(\int_{\mR^d}\Big(\cL^{\nu_{r,a(r,X_r)}}_{\sigma_{r,a(r,X_r)}}u(r,z)
-\cL^{\nu_{r,a(r,z)}}_{\sigma_{r,a(r,z)}}u(r,z)\Big)\varrho_\eps(X_r-z)\dif z\right)\\
&\quad+1_{\alpha=1}\mE\left(\int_{\mR^d} ( b_{r,X_r}- b_{r,z})\cdot\nabla u(s,z) \varrho_\eps(X_r-z)\dif z\right),
\end{align*}
and
$$
\Lambda'_\eps(r):=\mE\left(\int_{\mR^d}\cJ^{(\bar\alpha)}_{u_\eps(r)}(X_r,\bar\sigma_{r,X_r}y)\bar\nu_{r,X_r}(\dif y)\right)
+1_{\alpha\in(1,2)}\mE\left((\bar b\cdot\nabla u_\eps)(r,X_r)\right).
$$
Hence, by the dominated convergence theorem, we have
\begin{align}\label{ED1}
\begin{split}
&\mE\left(\int^T_0f(r,X_r)\dif r\right)=\lim_{\eps\to 0}\mE\left(\int^T_0f_\eps(r,X_r)\dif r\right)\\
&\quad\leq\sup_{\eps\in(0,1)}\left(\|u_\eps(0)\|_\infty+\int^T_0 \Big(|\Lambda_\eps(r)|+|\Lambda'_\eps(r)|\Big)\dif r\right).
\end{split}
\end{align}
Since $\int_{\mR^d}\varrho_\eps(z)\dif z=1$,
by Jensen's inequality and Lemma \ref{Le31}, recalling 
$q>\frac{d}{\gamma_0(\gamma_\sigma(\alpha\wedge 1)\wedge\gamma_\nu)}\vee\frac{m}{\gamma_\sigma(\alpha\wedge 1)\wedge\gamma_\nu}$, we have
\begin{align}
|\Lambda_\eps(r)|&\leq\mE\left(\int_{\mR^d}|\cL^{\nu_{r,a(r,X_r)}}_{\sigma_{r,a(r,X_r)}}u(r,z)
-\cL^{\nu_{r,a(r,z)}}_{\sigma_{r,a(r,z)}}u(r,z)|^q\varrho_\eps(X_r-z)\dif z\right)^{1/q}\no\\
&\quad+1_{\alpha=1}\mE\left(\int_{\mR^d}|( b_{r,X_r}- b_{r,z})\cdot\nabla u(r,z)|^q\varrho_\eps(x-z)\dif z\right)^{1/q}\no\\
&\leq\sup_x\left(\int_{\mR^d}|\cL^{\nu_{r,a(r,x)}}_{\sigma_{r,a(r,x)}}u(r,z)-\cL^{\nu_{r,a(r,z)}}_{\sigma_{r,a(r,z)}}u(r,z)|^q\varrho_\eps(x-z)\dif z\right)^{1/q}\no\\
&\quad+1_{\alpha=1}\sup_x\left(\int_{\mR^d}|( b_{r,x}- b_{r,z})\cdot\nabla u(r,z)|^q\varrho_\eps(x-z)\dif z\right)^{1/q}\no\\
&\leq\|\varrho_\eps\|^{1/q}_\infty\Big(C\hbar_a(\eps)^{\gamma_\sigma\beta\wedge\gamma_\nu}
\|\Delta^{\frac{\alpha}{2}} u(r)\|_q+ 1_{\alpha=1}\hbar_b(\eps)\|\nabla u(r)\|_q\Big)\no\\
&\leq C\eps^{\gamma_0(\gamma_\sigma\beta\wedge\gamma_\nu)-d/q}\|\Delta^{\frac{\alpha}{2}} u(r)\|_q
+C1_{\alpha=1}\eps^{\gamma_0-d/q}\|\Delta^{\frac{1}{2}} u(r)\|_q,\label{ED2}
\end{align}
where $\beta\in(\frac{m}{p\gamma_\sigma}\vee\frac{d}{q\gamma_0\gamma_\sigma},\alpha\wedge1)$.
Similarly, recalling \eqref{AS2}, by Lemma \ref{Le41}, \eqref{Es304} and Lemma \ref{Le23}, we have
\begin{align}
|\Lambda'_\eps(r)|&\leq 
\sup_x\left|\int_{\mR^d}\cJ^{(\bar\alpha)}_{u_\eps(r)}(x,\bar\sigma_{r,x}y)\bar\nu_{r,x}(\dif y)\right|
+1_{\alpha\in(1,2)}\|\bar b\|_\infty\|\nabla u_\eps(r)\|_\infty\no\\
&\preceq\|u_\eps(r)\|_{\theta,q}+1_{\alpha\in(1,2)}\|\bar b\|_\infty\|\nabla u(r)\|_\infty\preceq \|u(r)\|_{\theta,q}\no\\
&=\|u(r)-u(T)\|_{\theta,q}\preceq \|u\|_{\mU^{\alpha,q}(T)}\preceq \|f\|_{\mL^{q}(T)},\label{ED4}
\end{align}
Combining \eqref{ED1}-\eqref{ED4}, we obtain
$$
\mE\left(\int^T_0f(r,X_r)\dif r\right)\leq C\|f\|_{\mL^q(T)}+C\int^T_0\|\Delta^{\frac{\alpha}{2}} u(r)\|_q\dif r\leq C\|f\|_{\mL^q(T)}.
$$
Thus, we obtain (\ref{EE89}).
\\
\\
(2) Let $0\leq t_1<t_2\leq T$. For any $f\in \mL^p(T)\cap\mL^q(T)$, let $u\in\mU^{\alpha,p}(t_2)\cap \mU^{\alpha,q}(t_2)$ 
solve equation 
$$
\p_t u+\sL u=f,\ \  u(t_2)=0.
$$ 
Let $u_\eps(t,x):=u(t)*\varrho_\eps(x)$. Then 
$$
\p_t u_\eps+\sL u_\eps=f_\eps+\sL u_\eps-(\sL u)*\varrho_\eps,\ \ u_\eps(t_2)=0.
$$
By Remark \ref{Re} again, we have
\begin{align}
-u_\eps(t_1,X_{t_1})&=\mE\left(\int^{t_2}_{t_1}(\p_r u_\eps(r,X_r)+\sL_r u_\eps(r,X_r))\dif r\Big|\sF_{t_1}\right)\no\\
&=\mE\left(\int^{t_2}_{t_1}f_\eps(r,X_r)\dif r\Big|\sF_{t_1}\right)+\xi^\eps_{t_1,t_2},\label{EU1}
\end{align}
where
$$
\xi^\eps_{t_1,t_2}:=\mE\left(\int^{t_2}_{t_1}\big(\sL u_\eps-(\sL u)*\varrho_\eps\big)(r,X_r)\dif r\Big|\sF_{t_1}\right).
$$
By (\ref{EE89}), we have
$$
\lim_{\eps\to 0}\mE\left(\int^{t_2}_{t_1}|f_\eps(r,X_r)-f(r,X_r)|\dif r\right)\leq C\lim_{\eps\to 0}\|f_\eps-f\|_{\mL^q(T)}=0,
$$
and by (\ref{Es17}) and (\ref{Es10}),
$$
\lim_{\eps\to 0}\mE|\xi^\eps_{t_1,t_2}|\leq C\lim_{\eps\to 0}\|\sL u_\eps-(\sL u)*\varrho_\eps\|_{\mL^q(T)}=0.
$$
By taking limits for both sides of (\ref{EU1}), we get
\begin{align}\label{NG1}
-u(t_1,X_{t_1})=\mE\left(\int^{t_2}_{t_1}f(r,X_r)\dif r\Big|\sF_{t_1}\right).
\end{align}
By Lemma \ref{Le23} with $\beta\in(\frac{d}{p},\alpha(1-\frac{1}{p}))$ and \eqref{Es88}, we have
\begin{align*}
\|u(t_1)\|_\infty&\leq C\|u(t_1)\|_{\beta,p}=C\|u(t_1)-u(t_2)\|_{\beta,p}\\
&\leq C(t_2-t_1)^{1-\frac{\beta}{\alpha}-\frac{1}{p}}\|u\|_{\mU^{\alpha,p}(T)}\leq C(t_2-t_1)^{1-\frac{\beta}{\alpha}-\frac{1}{p}}\|f\|_{\mL^{p}(s,T)}.
\end{align*}
Substituting this into \eqref{NG1}, we obtain (\ref{Eq24}).
\end{proof}
\bl\label{Le44}
Under {\bf (H$_A$)} and {\bf (H$'_B$)}, for each $(s,x)\in\mR_+\times\mR^d$, the set $\tilde \Gamma^\sL_{s,x}$ has at most one element.
\el
\begin{proof}
Let $\mP_1,\mP_2\in \Gamma^\sL_{s,x}$ satisfy that for some $p_1>1$ and any $p\geq p_1$ and $f\in\mL^p(T)$,
$$
\mE^{\mP_i}\left(\int^T_s f(r,X_r)\dif r\right)\leq C\|f\|_{\mL^p(T)},\ \ i=1,2.
$$ 
Let $T>s$ and $p>\tfrac{m}{\gamma_\sigma(\alpha\wedge 1)\wedge\gamma_\nu}\vee\frac{d}{\alpha\wedge 1}\vee(\frac{d}{\alpha}+1)\vee\frac{d}{\bar\alpha}\vee p_1$. 
For any $t\in[s,T]$ and $f\in C_c(\mR^d)$, let $u\in\mU^{\alpha,p}(t)$ solve equation 
$$
\p_r u+\sL_r u=f,\ \  u(t)=0.
$$ 
As in the proof of step (3) of Theorem \ref{Main},  we have
$$
u(s,x)=\mE^{\mP_i}\left(\int^t_s f(r,X_r)\dif r\right),\  i=1,2,
$$
which implies that for any $f\in C_c(\mR^d)$,
$$
\mE^{\mP_1}f(X_t)=\mE^{\mP_2}f(X_t).
$$
In particular, for any $E\in\sB(\mR^d)$ and $t\geq s$,
\begin{align}\label{EE0}
\mP_1(X_t\in E)=\mP_2(X_t\in E).
\end{align}
Now let $\mP_1,\mP_2\in \tilde \Gamma^\sL_{s,x}$. Below we use induction to show that for any $s\leq t_1<t_2<\cdots<t_n\leq T$ and $E_1,\cdots,E_n\in\sB(\mR^d)$,
\begin{align}\label{EE1}
\mP_1(X_{t_1}\in E_1,\cdots,X_{t_n}\in E_n)=\mP_2(X_{t_1}\in E_1,\cdots,X_{t_n}\in E_n).
\end{align}
For $n=1$, it has been proven in \eqref{EE0}. Suppose that it holds for some $n$. By \cite[Theorem 1.2]{St}, the regular conditional probabilities
$\mP^\omega_1, \mP^\omega_2$ of $\mP_1$ and $\mP_2$ with respect to $\sG_{n}:=\sigma(X_{t_1},\cdots, X_{t_n})\subset\sF_{t_n}$ 
belong to $\Gamma^{\sL}_{t_n, X_{t_n}(\omega)}$ and satisfy
that for some $p_1>1$ and any $p\geq p_1$ and $f\in\mL^p(T)$,
$$
\mE^{\mP^\omega_i}\left(\int^T_{t_n} f(r,X_r)\dif r\right)
=\mE^{\mP_i}\left(\int^T_{t_n} f(r,X_r)\dif r\Big|\sG_{n}\right)(\omega)\leq C\|f\|_{\mL^p(T)}\ a.s.
$$ 
Notice that
\begin{align*}
&\mP_i(X_{t_1}\in E_1,\cdots,X_{t_n}\in E_n, X_{t_{n+1}}\in E_{n+1})\\
&\quad=\mE^{\mP_i}\left(1_{X_{t_1}\in E_1,\cdots,X_{t_n}\in E_n}\mP^\omega_i(X_{t_{n+1}}\in E_{n+1})\right),\ i=1,2.
\end{align*}
By the induction hypothesis and \eqref{EE0}, we get \eqref{EE1} for $n+1$.
\end{proof}
To show the existence of a martingale solution, we shall use the weak convergence argument. 
Let $\varrho:\mR^d\to\mR_+$ be a smooth function with support in the unit ball and $\int\varrho=1$.
For $n\in\mN$, let $\varrho_n(x)=n^d\varrho(nx)$ and define
$$
a_n(t,x):=a(t,\cdot)*\varrho_n(x),\ \  b^n_{t,x}:= b_{t,\cdot}*\varrho_n(x).
$$
We also assume that
\begin{enumerate}[{\bf (H$'_A$)}]
\item It holds that for some $\gamma_0$,
$$
\sup_{\eps\in(0,1)}(\hbar_{a_n}(\eps)+\hbar_{b_n}(\eps))\eps^{-\gamma_0}<\infty
$$
and
$$
\sup_{n,t}\|\sigma_{t,a_n(t,\cdot)}\|_\infty<\infty,\ \inf_{n,t}\inf_{x\in\mR^d}\inf_{|\xi|=1}|\sigma_{t,a_n(t,x)}\xi|>0.
$$
\end{enumerate}
Notice that under {\bf (H$_A$)}, {\bf (H$'_A$)} is automatically satisfied for the examples in ({\bf Ex}).
\bt
Assume that {\bf (H$_A$)}, {\bf (H$'_A$)}  and {\bf (H$'_B$)} hold. Then for any $(s,x)\in\mR_+\times\mR^d$,
$\tilde\Gamma^\sL_{s,x}$ has one and only one element.
\et
\begin{proof}
By Lemma \ref{Le44}, it suffices to show the existence. Define
$$
\bar\sigma^n_{t,x}:=\bar\sigma_{t,\cdot}*\varrho_n(x),\ \ \bar\nu^n_{t,x}:=\bar\nu_{t,\cdot}*\varrho_n(x),\ \ \bar b^n_{t,x}:=\bar b_{t,\cdot}*\varrho_n(x).
$$
It is easy to see that {\bf (H$'_B$)} holds uniformly with respect to $n$.
Let $A^n_{t,x}$ and $B^n_{t,x}$ be defined in terms of $a^n, b^n$ and $\bar\sigma^n,\bar\nu^n,\bar b^n$ respectively. Let $\sL^n:=A^n+B^n$.
By \cite[Chapter IX, Theorem 2.31]{Ja-Sh}, for each $(s,x)\in\mR_+\times\mR^d$, 
there exists at least one solution $\mP_n\in\Gamma^{\sL^n}_{s,x}$. 
Let $T>0$ and
$$
p>\tfrac{m}{\gamma_\sigma(\alpha\wedge 1)\wedge\gamma_\nu}\vee\tfrac{d}{\bar\alpha}\vee(\tfrac{d}{\alpha}+1)\vee\tfrac{d}{\alpha\wedge 1}.
$$
By Theorem \ref{Main}, there exists a constant $C$ independent of $n$ such that for all $s\leq t_1<t_2\leq T$ and $f\in\mL^p(T)$,
\begin{align}\label{NG3}
\mE^{\mP_n}\left(\int^{t_2}_{t_1}f(r,X_r)\dif r\Big|\sF_{t_1}\right)\leq C\|f\|_{\mL^p(s,T)}.
\end{align}
By \cite[Chapter IV, Theorem 4.18]{Ja-Sh}, $(\mP_n)_{n\in\mN}$ is tight. Let $\mP$ be an accumulation point of $(\mP_n)_{n\in\mN}$. We want to show that $\mP\in\tilde\Gamma^{\sL}_{s,x}$. Up to extracting a subsequence, we may assume that $\mP_n$ weakly converges to $\mP$.
For $f\in C_c((0,T)\times\mR^d)$, by taking weak limits for \eqref{NG3}, we have
$$
\mE^{\mP}\left(\int^{t_2}_{t_1}f(r,X_r)\dif r\Big|\sF_{t_1}\right)\leq C\|f\|_{\mL^p(s,T)}.
$$
By a standard monotone class argument, the above estimate still holds for all $f\in\mL^p(T)$.

It remains to show $\mP\in\Gamma^{\sL}_{s,x}$. Let $\varphi\in C^\infty_c(\mR^d)$. It suffices to show that
for any $s\leq t_1\leq t_2\leq T$ and any bounded $\sF_{t_1}$-measurable continuous functional $G$, 
$$
\mE^{\mP}\big(G M^\varphi_{t_2}\big)=\mE^{\mP}\big(G M^\varphi_{t_1}\big).
$$
Since $\mP_n\in\Gamma^{\sL^n}_{s,x}$, we have
\begin{align}\label{NG2}
\mE^{\mP_n}\big(G M^{n,\varphi}_{t_2}\big)=\mE^{\mP_n}\big(G M^{n,\varphi}_{t_1}\big),
\end{align}
where $M^{n,\varphi}_{t_i}:=\varphi(X_{t_i})-\int^{t_i}_s\sL^n_{r}\varphi(X_r)\dif r$.
We naturally want to take limits for both sides of \eqref{NG2}, that is, to prove
\begin{align}
\lim_{n\to\infty}\mE^{\mP_n}\left(G \int^{t_i}_sA^n_r\varphi(X_r)\dif r\right)&=\mE^{\mP}\left(G \int^{t_i}_sA_r\varphi(X_r)\dif r\right),\label{NG5}\\
\lim_{n\to\infty}\mE^{\mP_n}\left(G \int^{t_i}_sB^n_r\varphi(X_r)\dif r\right)&=\mE^{\mP}\left(G \int^{t_i}_sB_r\varphi(X_r)\dif r\right).\label{NG7}
\end{align}
Let us only prove \eqref{NG7} since \eqref{NG5} is similar.
For each $m\in\mN$,  since $\varphi\in C^\infty_c(\mR^d)$, by definition \eqref{AB2}, it is easy to see that the functional 
$$
\omega\mapsto \int^{t_i}_sB^m_r\varphi(X_r(\omega))\dif r=\int^{t_i}_sB^m_r\varphi(\omega_r)\dif r
$$
is bounded and continuous. Thus,
\begin{align}\label{NG4}
\lim_{n\to\infty}\mE^{\mP_n}\left(G \int^{t_i}_sB^m_r\varphi(X_r)\dif r\right)=\mE^{\mP}\left(G \int^{t_i}_sB^m_r\varphi(X_r)\dif r\right).
\end{align}
On the other hand, by \eqref{NG3}, Lemmas \ref{Le41}, \ref{Le22} and the dominated convergence theorem, we have
\begin{align*}
&\lim_{m\to\infty}\sup_n\left|\mE^{\mP_n}\left(G \int^{t_i}_s(B^m_r\varphi-B_r\varphi)(X_r)\dif r\right)\right|^p\\
&\preceq \|G\|_\infty\lim_{m\to\infty}\int^{t_i}_s\!\!\!\int_{\mR^d}|B^m_r\varphi(x)-B_r\varphi(x)|^p\dif x\dif r=0,
\end{align*}
which together with \eqref{NG4} implies \eqref{NG7}. The proof is complete.
 \end{proof}
 
 \bc
 Consider SDE \eqref{SDE0}. Suppose that $\sigma(x)$ is bounded and uniformly continuous and nondegenerate, $\bar\sigma(x)$ is bounded measurable,
 $L_\cdot$ is an $\alpha$-stable L\'evy process with L\'evy measure $\nu\in\bL^{(\alpha)}_{non}$, and  
 $\bar L_\cdot$ is a $\beta$-stable L\'evy process with $\beta<\alpha$, and independent of $L_\cdot$. For any $x\in\mR^d$, 
 there exists a unique weak solution $X_t$ to \eqref{SDE0}
 with the property that for any $p>\frac{d^2}{\alpha\wedge 1}\vee(\frac{d}{\alpha}+1)\vee\frac{d}{\beta}$ and $T>0$, and all $0\leq t_1\leq t_2\leq T$ and $f\in\mL^p(T)$,
\begin{align}\label{Prop}
\mE\left(\int^{t_2}_{t_1} f(s,X_s)\dif s\Big|\sF_{t_1}\right)\leq C\|f\|_{\mL^p(T)}.
\end{align}
 \ec
 \begin{proof}
The existence of a weak solution with property \eqref{Prop} follows by Theorem \ref{Main} and
the same argument as in \cite{Ba-Ch}. Since the law of any weak solution of SDE \eqref{SDE0} 
with property \eqref{Prop} belongs to $\tilde\Gamma^{\sL}_{0,x}$, the uniqueness follows by Lemma \ref{Le44}.
 \end{proof}
 
\section{Pathwise uniqueness of SDEs driven by L\'evy processes}

In this section we prove a pathwise uniqueness result for SDE (\ref{SDE}). We recall the following simple result (cf. \cite[Lemma 2.6]{Zh1}).
\bl\label{Le0}
Let $(Z_t)_{t\geq 0}$ be a locally bounded and ($\sF_t$)-adapted process and $(\ell_t)_{t\geq 0}$
a continuous real valued non-decreasing ($\sF_t$)-adapted process with $\ell_0=0$.
Assume that for any stopping time $\tau$ and $t\geq 0$,
$$
\mE|Z_{t\wedge\tau}|\leq\mE\int^{t\wedge\tau}_0|Z_s|\dif \ell_s.
$$
Then $Z_t=0$ a.s. for all $t\geq 0$.
\el
We also need the following elementary inequality.
\bl\label{Le52}
For any $q\in(0,1)$, there exists a constant $C=C(q)>0$ such that for all $x,y\in\mR^d$,
\begin{align}
|x|x|^{q-1}-y|y|^{q-1}|\leq C|x-y|^q.\label{EP}
\end{align}
\el
\begin{proof}
Let $\bar y=y/|y|$. It suffices to prove
$$
|x|x|^{q-1}-\bar y|\leq C|x-\bar y|^q,\ \ \forall x,y\in\mR^d.
$$
By the coordinate rotation, we can assume $\bar y=(1,0,\cdots,0)$ so that it suffices to prove that for any $x=(x_1,\cdots,x_d)$,
$$
|x_1|x|^{q-1}-1|^2+\sum_{i=2}^d|x_i|^2|x|^{2(q-1)}\leq C^2 \left(|x_1-1|^2+\sum_{i=2}^d |x_i|^2\right)^q,
$$
which is equivalent to prove
$$
|x|^{2q}-2x_1|x|^{q-1}+1\leq C^2 (|x|^2+1-2x_1)^q.
$$
Define
$$
f_a(b):=C^2(a^2+1-2b)^q-a^{2q}+2ba^{q-1}-1,\ \ |b|\leq a.
$$
Clearly, if $C\geq 2$, then
\begin{align}
f_a(a)\geq 0,\ \ f_a(-a)\geq 0.\label{EA1}
\end{align}
Now we consider the minimal point of $b\mapsto f_a(b)$ on $[-a,a]$. Solving the following equation,
$$
0=f'_a(b)=-2C^2q (a^2+1-2b)^{q-1}+2a^{q-1},
$$
we obtain
$$
b_0=\frac{a^2+1-2\beta a}{2},\ \ \beta:=(C^2q)^{\frac{1}{1-q}}/2.
$$
Since $|b_0|\leq a\Rightarrow |a^2+1-2\beta a|\leq 2a$, we have
$$
\gamma_1:=\beta+1-\sqrt{\beta^2+2\beta}\leq a\leq\beta-1-\sqrt{\beta^2-2\beta}=:\gamma_2,
$$ 
or
$$
\gamma_3:=\beta-1+\sqrt{\beta^2-2\beta}\leq a\leq \beta+1+\sqrt{\beta^2+2\beta}=:\gamma_4.
$$
If $0\leq a\notin[\gamma_1,\gamma_2]\cup[\gamma_3,\gamma_4]$, then there is no zero points for $f'_a(b)=0$, and by (\ref{EA1}) we have
$$
f_a(b)\geq 0,\ \ \forall b\in[-a,a].
$$
It remains to prove that if $C$ is large enough, then for any $a\in[\gamma_1,\gamma_2]\cup[\gamma_3,\gamma_4]$, 
$$
f_a(b_0)\geq 0.
$$
Notice that
$$
f_a(b_0)=a^{q+1}-a^{2q}+C_q a^q+a^{q-1}-1,
$$
where
$$
C_q:=C^2(C^2q)^{q/(1-q)}-(C^2q)^{1/(1-q)}
=C^{2/(1-q)}q^{q/(1-q)}(1-q)>0.
$$
Since $C\mapsto\beta(C)$ is increasing and $q\in(0,1)$, one sees that if $C$ is large enough, then for any $a\in[\gamma_3,\gamma_4]$, 
$$
f_a(b_0)\geq a^{q+1}-a^{2q}-1\geq 0;
$$
while for $a\in[\gamma_1,\gamma_2]$, since $\gamma_2(\beta)\leq\frac{1}{2\sqrt{\beta^2-2\beta}}\to 0$ as $\beta\to\infty$, we also have
$$
f_a(b_0)\geq a^{q-1}-a^{2q}-1\geq 0.
$$
The desired inequality follows.
\end{proof}
\br
If the constant $C$ in (\ref{EP}) is allowed to be dependent on the dimension $d$, then we have the following simple proof:
Without loss of generality, we assume $|y|\leq|x|$. First of all, we assume $|x-y|>\frac{|x|}{2}$. In this case, we have
$$
|x|x|^{q-1}-y|y|^{q-1}|\leq |x|^q+|y|^q\leq 2|x|^q\leq 2^{1+q}|x-y|^q.
$$
Next, we assume $|x-y|\leq\frac{|x|}{2}$. If we set $f(x)=x|x|^{q-1}$, then
\begin{align*}
|f(x)-f(y)|&\leq|x-y|\int^1_0|\nabla f(x+\theta(y-x))|\dif \theta\\
&\leq C_{d,q}|x-y|\int^1_0|x+\theta(y-x)|^{q-1}\dif \theta\\
&\leq C_{d,q}|x-y|(|x|-|y-x|)^{q-1}\leq C_{d,q}|x-y|^q.
\end{align*}
\er
Now we can prove the following main result of this section.
\bt
Let $L_t$ be a symmetric L\'evy process with  L\'evy measure $\nu$.
Suppose that for some $\alpha\in(0,2)$ and $\nu^{(\alpha)}_1, \nu^{(\alpha)}_2\in\bL^{(\alpha)}_{non}$,
$$
\nu^{(\alpha)}_1\leq \nu\leq \nu^{(\alpha)}_2,
$$ 
and $\sigma(x)$ is linear growth and nondegenerate, and for some $p>d(1+\alpha\vee 1)$,
$$
\sigma(x)\in\mW^{1,p}_{loc}.
$$
Then for any $x\in\mR^d$, there exists a unique strong solution to the following SDE:
\begin{align}
\dif X_t=\sigma(X_{t-})\dif L_t,\ \ X_0=x.\label{SDE1}
\end{align} 
\et
\begin{proof}
By a standard localization argument, we can assume that
$$
\mbox{$\sigma\in\mW^{1,p}$ is bounded and uniformly nondegenerate}.
$$
In this case, since $\sigma$ is continuous, the existence of a weak solution is standard. 
We only show the pathwise uniqueness.
Let $N(\dif t,\dif z)$ be the Poisson random measure associated with $L_t$, i.e.,
$$
N(t,\Gamma):=\sum_{s\leq t}1_{\Gamma}(\Delta L_s).
$$
By L\'evy-It\^o's decomposition, we have
$$
L_t=\int_{|z|\leq 1}y\tilde N(t,\dif z)+\int_{|y|>1}yN(t,\dif z).
$$
Thus, SDE (\ref{SDE1}) can be written as
$$
X_t=x+\int^t_0\!\!\!\int_{|z|\leq 1}\sigma(X_{s-})z\tilde N(\dif s,\dif z)+\int^t_0\!\!\!\int_{|z|>1}\sigma(X_{s-})z N(\dif s,\dif z).
$$
Let $X_t$ and $Y_t$ be two solutions of SDE (\ref{SDE1}) starting from the same point $x\in\mR^d$. Set 
$$
Z_t:=X_t-Y_t.
$$
For $\eps\geq 0$ and $q>0$, let $f_\eps(x):=(|x|^2+\eps)^{q/2}$. By It\^o's formula, we have
\begin{align}
f_\eps(Z_t)&=\int^t_0\!\!\!\int_{|z|\leq 1}[f_\eps(Z_{s-}+\Sigma_s z)-f_\eps(Z_{s-})-\Sigma_s z\cdot\nabla f_\eps(X_{s-})]\nu(\dif z)\dif s\no\\
&\quad+\int^t_0\!\!\!\int_{|z|\leq 1}[f_\eps(Z_{s-}+\Sigma_s z)-f_\eps(Z_{s-})]\tilde N(\dif s,\dif z)\no\\
&\quad+\int^t_0\!\!\!\int_{|z|>1}[f_\eps(Z_{s-}+\Sigma_s z)-f_\eps(Z_{s-})]N(\dif s,\dif z),\label{NB1}
\end{align}
where
$$
\Sigma_s:=\sigma(X_{s-})-\sigma(Y_{s-}).
$$
Let $\tau_0=0$ and define recursively for $n\in\mN$,
$$
\tau_n:=\inf\{s>\tau_{n-1}: |\Delta L_s|\geq 1\}.
$$
Let $\tau$ be any bounded stopping time. 
By (\ref{NB1}), we have
\begin{align*}
\mE f_\eps(Z_{t\wedge\tau_1\wedge\tau})&=\mE\int^{t\wedge\tau_1\wedge\tau}_0\!\!\!\!\int_{|z|\leq 1}\!\!
[f_\eps(Z_{s-}+\Sigma_s z)-f_\eps(Z_{s-})-\Sigma_s z\cdot\nabla f_\eps(X_{s-})]\nu(\dif z)\dif s.
\end{align*}
(Case: $\alpha\in(0,1)$) Let $q\in(\alpha,1)$. By the symmetry of $\nu$ and letting $\eps\to 0$, we have
\begin{align}
\mE |Z_{t\wedge\tau_1\wedge\tau}|^q&=\mE\int^{t\wedge\tau_1\wedge\tau}_0\!\!\!\int_{|z|\leq 1}[|Z_{s-}+\Sigma_s z|^q-|Z_{s-}|^q]\nu(\dif z)\dif s\no\\
&\leq\mE\int^{t\wedge\tau_1\wedge\tau}_0\!\!\!\int_{|z|\leq 1}|\Sigma_s z|^q\nu(\dif z)\dif s\quad (\because ||x|^q-|y|^q|\leq|x-y|^q)\no\\
&\leq C\mE\left(\int^{t\wedge\tau_1\wedge\tau}_0|Z_s|^q\dif \ell_s\right)\int_{|z|\leq 1}|z|^q\nu^{(\alpha)}_2(\dif z),\label{EA2}
\end{align}
where the last step is due to Lemma \ref{Le21} and
$$
\ell_t:=\int^t_0(\cM|\nabla\sigma|(X_{s-})+\cM|\nabla\sigma|(Y_{s-}))^q\dif s.
$$
(Case: $\alpha\in[1,2)$) Let $q\in(\alpha,2)$.  By Lemma \ref{Le52}, we have
\begin{align}
\mE f_0(Z_{t\wedge\tau_1\wedge\tau})&=\mE\int^{t\wedge\tau_1\wedge\tau}_0\!\!\!\!\int_{|z|\leq 1}\!\Sigma_s z\cdot
\!\int^1_0[\nabla f_0(Z_{s-}+\theta\Sigma_s z)-\nabla f_0(Z_{s-})]\dif\theta\nu(\dif z)\dif s\no\\
&\leq C\mE\int^{t\wedge\tau_1\wedge\tau}_0\!\!\!\int_{|z|\leq 1}|\Sigma_s z|^q\nu(\dif z)\dif s\int^1_0\theta^{q-1}\dif \theta\no\\
&\leq C\mE\left(\int^{t\wedge\tau_1\wedge\tau}_0|Z_s|^q\dif A_s\right)\int_{|z|\leq 1}|z|^q\nu^{(\alpha)}_2(\dif z).\label{EA3}
\end{align}
Since $\sigma\in\mW^{1,p}$ with $p>d(1+\alpha\vee 1)$, by Sobolev's embedding theorem, we have
$$
|\sigma(x)-\sigma(y)|\leq C|x-y|^\gamma,\ \ \forall\gamma\in(0,1-\tfrac{d}{p}).
$$
If we choose $\gamma$ close to $1-\frac{d}{p}$ and $q$ close to $\alpha$ so that
$$
\tfrac{p}{q}>\tfrac{d}{\gamma (\alpha\wedge 1)}\vee(\tfrac{d}{\alpha}+1),
$$
then by (\ref{Eq24}) and Lemma \ref{Le21}, we have
$$
\mE \ell_t\leq C\|(\cM|\nabla\sigma|)^q\|_{p/q}=C\|\cM|\nabla\sigma|\|^q_{p}\leq C\|\nabla\sigma\|^q_{p}.
$$
By (\ref{EA2}), (\ref{EA3}) and Lemma \ref{Le0}, we obtain $Z_{t\wedge\tau_1}=0$. 
Proceeding the above proof, we can prove $Z_{t\wedge\tau_n}=0$ for any $n$. The proof is complete.
\end{proof}

\end{document}